\documentclass[12pt,oneside, reqno]{amsart}
 \usepackage{graphicx}
 \usepackage{hhline}
\usepackage{amsmath,amsthm,amsfonts,amscd,amssymb,mathrsfs}
\usepackage{xspace}
\usepackage[all]{xypic}
\usepackage{booktabs} 
\usepackage{physics}
\usepackage[numeric,backrefs,lite,initials]{amsrefs}
\usepackage[backref]{hyperref}
\usepackage{soul} 

\usepackage{verbatim}
\usepackage{amscd}   
\usepackage[all]{xy} 
\usepackage{youngtab} 
\usepackage{ytableau} 

\usepackage{nicefrac}
\usepackage{xfrac}

\usepackage{tikz}

\newcommand{\ccircle}[1]{* + <1ex>[o][F-]{#1}}
\newcommand{\ccirc}[1]{\xymatrix@1{* + <1ex>[o][F-]{#1}}}

\usepackage[T1]{fontenc} 
\usepackage{cleveref} 
\numberwithin{equation}{section} 

\usepackage{color}

\makeatletter
\newtheorem{rep@theorem}{\rep@title}
\newcommand{\newreptheorem}[2]{%
\newenvironment{rep#1}[1]{%
 \def\rep@title{#2 \ref{##1}}%
 \begin{rep@theorem}}%
 {\end{rep@theorem}}}
\makeatother

\newtheorem{theorem}{Theorem}[section]
\newreptheorem{theorem}{Theorem}
\newreptheorem{lemma}{Lemma}

\newtheorem{conjecture}[theorem]{Conjecture}

\newcommand{\SL}{\operatorname{SL}}
\newcommand{\GL}{\operatorname{GL}}

\newcommand{\diag}{\operatorname{diag}}

\newcommand{\CC}{\mathbb{C}}
\newcommand{\ZZ}{\mathbb{Z}}

\def\bw#1{{\textstyle\bigwedge^{\hspace{-.2em}#1}}}

\def\phi{ \varphi }
\def\ep{ \varepsilon}

\def \a{\alpha}
\def \b{\beta}
\def \n{\mathfrak{n}}
\def \h{\mathfrak{h}}
\def \d{\mathfrak{d}}
\def \z{\mathfrak{z}}

\def \c{\mathfrak{c}}
\def \s{\mathfrak{s}}
\def \ga{\gamma}

\def \g{\mathfrak{g}}

\def \sl{\mathfrak{sl}}
\def \diag{\textrm{diag}}

\theoremstyle{definition}

\newtheorem{examplex}{Example}
\newenvironment{example}
  {\pushQED{\qed}\examplex}
  {\popQED\endexamplex}

\theoremstyle{remark}

\newcommand{\ad}{\operatorname{ad }}

\newcommand{\ch}[7]{#1\,#2\,#3\,#4\,#5\,#6\,#7}

\begin{document}
\date{\today}

\author{Luke Oeding}\email{oeding@auburn.edu}
\address{Department of Mathematics and Statistics,
Auburn University,
Auburn, AL, USA
}

\title[Antonyan's classification of $\SL_8$ orbits in $\bw{4}\CC^8$]{A Translation of ``Classification of four-vectors of an 8-dimensional space,'' by Antonyan, L. V. , with an appendix by the translator}
\begin{abstract} 
We give a translation of the entitled paper \cite{Antonyan}.  We include an appendix that shows how to produce normal forms for each nilpotent orbit. 
\end{abstract}
\maketitle

\section*{Author's Introduction (provided by AMS's MathSciNet Entry)}
By classification of four-vectors of an eight-dimensional space we mean a description of the orbits of a natural representation in the space of four-vectors of a group of unimodular linear transformations of a fundamental eight-dimensional space. An analogous problem was solved by \`E. B. Vinberg and A. G. {\`E}la{\v{s}}vili, in a paper \cite{Vinberg-Elasvili} on the classification of trivectors of a nine-dimensional space. In the present paper we use the method of classifying nilpotent elements of graded Lie algebras proposed by them and perfected by Vinberg in another paper \cite{Vinberg76}.

For classification of four-vectors of an eight-dimensional space we use their realization in a $\mathbb{Z}_2$-graded simple Lie algebra of type $E_7$ and the general results of Vinberg in a third paper \cite{Vinberg75} on orbits and invariants of linear groups associated with graded Lie algebras.

\section{Basic Definitions}
\subsection{Some Notation}
We will use the following notation: If $V$ is a vector space, then $V^*$ is the dual space, $\bigotimes^k V$, $\bw{k}V$, $S^k V$ are respectively the $k$-th tensor, exterior, and symmetric powers of the vector space $V$,  $L(V) = V\otimes V^*$ is the space of linear transformations of the space $V$, $L_0(V)$ is the subspace of  $L(V)$ consisting of the transformations with trace $0$.

We will identify the space $\bw{k}V$ with a subspace in $\bigotimes^k V$, assuming that for $v_1,\ldots, v_k \in V$ 
\begin{equation}\label{eq:wedge}
v_1 \wedge \ldots \wedge v_k = \sum_\sigma (-1)^\sigma v_{\sigma_1} \otimes \cdots \otimes v_{\sigma_k}
\end{equation}
(where the sum on the right side of \eqref{eq:wedge} is taken over all permutations $\sigma = (\sigma_1,\ldots,\sigma_k)$ of the numbers $1,\ldots, k$, and $(-1)^\sigma$ is the parity of the permutation $\sigma$).

\subsection{Four-vectors of an eight-dimensional space}\label{sec:an12} 
From now on, by the letter $V$ we will denote the eight-dimensional vector space over the field of complex numbers. 

The four-vectors of the space $V$ are elements of the space $\bw{4}V$. Our task is to describe the orbits of the natural representation in the space of the four-vectors of the group of unimodular linear transformations of the base space $V$. In accordance with this, four-vectors $u_1$, $u_2$ are called equivalent if there exists a transformation $a\in \SL(V)$ so that $u_2=a(u_1)$. 

The stabilizer subgroup of four-vectors is called the subgroup 
\begin{equation}
S(u)=\{a \in \SL(V)\mid a(u)=u\}\subset \SL(V), 
\end{equation}
and the associated Lie algebra is the stabilizer subalgebra.

\subsection{The Algebra $E_7$}\label{sec:e7}
Let us introduce on a vector space 
\begin{equation}\label{eq:graded}
\mathfrak{g} = L_0(V) \oplus \bw{4} V
\end{equation}
the structure of a Lie algebra that is graded modulo 2 in such a way that the direct summands in \eqref{eq:graded} are the graded subspaces. 
Namely, on $L_0(V)$ we introduce the usual commutator bracket on operators, which, in tensor notation, is given by the following formula: 
\begin{equation}\label{eq:LLaction}
 [A, B]^i_j=A^i_kB^k_j -B^i_kA^k_j  \quad (A, B \in L_0(V)).
\end{equation}
Let us define the brackets of elements from $L_0(V)$ with elements from $\bw{4}V$ in such a way that it defines the usual action of the algebra $L_0(V)$ in the space $\bw{4}V$, i.e., by the formula: 
\begin{equation}\label{eq:LWaction}
[A,T]^{ijkl} = A^i_s T^{sjkl} + A^j_s T^{iskl} + A^k_s T^{ijsl} + A^l_s T^{ijks} \quad (A \in L_0(V),\; T\in \bw{4}V).
\end{equation}
The bracket of elements of the space $\bw{4}V$ with each other is defined as follows: 
\begin{equation}\label{eq:WWaction}
[T_1,T_2]^k_m = -\frac{1}{288} \delta_{p_1 q_1 r_1 s_1 p_2 q_2 r_2 m} (T_1^{p_1q_1r_1s_1}T_2^{p_2 q_2 r_2 k} - T_2^{p_1q_1r_1s_1}T_1^{p_2 q_2 r_2 k} )
\quad (T_1, T_2 \in \bw{4} V),
\end{equation}
where $\delta$ is a fixed non-zero element of the space $\bw{8}V^*$. 

A direct verification shows that the space $\mathfrak{g}$ with the bracket given by formulas \cref{eq:LLaction,eq:LWaction,eq:WWaction} is indeed a graded Lie algebra, on which one can define a grading by setting $\mathfrak{g}_0= L_0(V)$, $\mathfrak{g}_1 = \bw{4}V$. 

Let us fix a basis $\{e_1,\ldots, e_8\}$ of the space $V$ satisfying the condition $\delta (e_1,\ldots, e_8) = 1$, and denote by $\mathfrak{d}$ the commutative
subalgebra of the algebra $L_0(V)$ consisting of linear transformations that when expressed with respect to the basis $\{e_1,\ldots, e_8\}$  are diagonal matrices. 
Let $\ep_1, \ldots, \ep_8$ denote linear functionals on $\mathfrak{d}$ defined by formulas:
\begin{equation}\label{eq:1.7}
\varepsilon_k(\diag(c_1,\ldots,c_8)) = c_k \quad (k=1,\ldots, 8)
\end{equation}
Obviously, $\sum_{k=1}^8 \ep_k = 0$. 

Let $\{e^1,\ldots, e^8\}$ be the basis of the space $V^*$ dual to the basis $\{e_1,\ldots, e_8\}$.
For $A\in \d$ the following relations hold: 
\begin{equation}
[A, e_i\otimes e^j] = (\ep_i - \ep_j) A \cdot e_i \otimes e^j,
\end{equation}
\begin{equation}\label{eq:bracket01}
[A, e_i\wedge e_j \wedge e_k \wedge e_l] = (\ep_i + \ep_j + \ep_k + \ep_l) A \cdot e_i\wedge e_j \wedge e_k \wedge e_l.
\end{equation}
Consequently, the one-dimensional subspaces spanned by the elements 
\[
e_i \otimes e^j,\quad e_i\wedge e_j \wedge e_k \wedge e_l
\]
($i,\; i,\; k,\; l$ pairwise distinct) are root subspaces with respect to $\d$. Direct calculation shows that the simple algebra $\g$ is of type $E_7$,  $\d$ is a Cartan subalgebra, and
\begin{equation}
\Sigma= \{ \ep_i-\ep_j, (\ep_i + \ep_j + \ep_k + \ep_l) \} \quad (i,\; i,\; k,\; l \;\;\text{pairwise distinct})
\end{equation}
is a root system relative to $\d$. The subalgebra $\g_0$ is a simple Lie algebra of type $A_7$, $\d$ is a Cartan subalgebra, and
$\Sigma_0 = \{\ep_i-\ep_j\}$   ($i,\; j$ distinct) is a root system with respect to $\d$.

\subsection{The automorphism $\theta$ and group $G_0$}
Consider the transformation $\theta \colon \g \to \g\theta$ 
 given by the formula: 
\begin{equation}
\theta(x) = (-1)^kx \quad \text{for} \;\; x \in \g_k,
\end{equation}
which is a second-order automorphism of the algebra $\g$, and 
\begin{equation}\label{eq:g0}
\g_0 = \{ x\in \g\mid \theta(x) = x\}. 
\end{equation}

Let $G$ be a connected group having $\g$ as its Lie algebra. The automorphism $\theta$ is induced by some uniquely defined automorphism $\Theta$ of the group $G$.  The subgroup $G_0$, consisting of all elements of the group $G$ that are fixed with respect to $\Theta$, is connected, and (by \eqref{eq:g0}) is locally isomorphic to the group $\SL(V)$. Hence, there exists an epimorphism 
$\pi \colon \SL(V) \to G_0$, which induces the identity mapping of the Lie algebras. 
The adjoint representation $\text{Ad}$ of the group $G$ induces, by restriction, a linear representation of the group $G_0$ on the space $\g$. The composition  $\text{Ad}\circ \pi$ is a linear representation on the space $\g$ of the group $SL(V)$, which coincides on $\g_0$ with the adjoint representation, and on $\g_1$ with the natural representation of $\SL(V)$ acting on the space $\bw{4}V$. 

\section{General Theory}\label{sec:2}
In this section, we will summarize some results contained in \cite{Vinberg76} and \cite{Vinberg75}. Let $\g$ be a complex semisimple Lie algebra graded modulo $m$:
\begin{equation} 
\g = \bigoplus_{k\in \ZZ_m} \g_k.
\end{equation}
$G$ is a connected Lie group that has the algebra $\g$ as its Lie algebra. $G_0$ is  a connected subgroup of $G$ corresponding to the subalgebra $\g_0$. Article \cite{Vinberg76} studies the orbits and invariants of the linear representation of the group $G_0$ on the space $\g_1$ induced (by restriction) by the adjoint representation of the group $G$. In \cite{Vinberg75}, a classification of nilpotent elements 
in the space $\g_1$ is given. (We note that the latter comprise a finite number of orbits). 

\subsection{Jordan decomposition}\label{sec:2.1} Every element of the space $\g_1$ can be uniquely represented as the sum of commuting semisimple and nilpotent elements contained in $\g_1$.
Such a decomposition is called the Jordan decomposition, and the corresponding components are the semisimple and nilpotent parts the given element. 
The following theorem holds. 
\begin{theorem}
An element $x \in \g_1$ is nilpotent (respectively, semisimple) if and only if the closure of its orbit by the action of the group $G$ contains 0 (respectively, when its orbit is closed).
\end{theorem}
If two elements of the space $\g_1$ are equivalent, then their semisimple parts are also equivalent: two elements with coinciding semisimple parts are equivalent if and only if their corresponding nilpotent parts are transformed one into the other by a transformation from the group $G_0$, leaving the semisimple part fixed.

 Let $p$ be an arbitrary semisimple element of the space $\g_1$. Its centralizer $\mathfrak{z}(p)$ in the algebra $\g$ is a reductive subalgebra of maximum rank: we denote by $\tilde \z (p)$ the commutator subalgebra of the subalgebra $\z(p)$. The grading of the algebra $\g$ induces the grading of the algebra $\z(p)$
 \[
 \tilde \z(p) = \bigoplus_{k\in\ZZ_m} \tilde \z_k(p),\quad \tilde \z_k(p) = \tilde \z (p) \cap \g_k.
 \]
 Every nilpotent element of the space $\g_1$  that commutes with $p$ lies in  $\tilde \z_1(p)$. 
 
 Let $Z(p)$ and $Z_0(p)$ be the centralizers (i.e., stabilizer subgroups for the conjugation action) of the element $p$ in the groups $G$ and $G_0$, respectively,  $\tilde Z_0(p)$ the projection of the subgroup $Z_0(p)$ onto the derived subgroup of the group $Z(p)$. In this notation, the classification (with respect to the group $G_0$) of the elements in the space $\g_1$ whose semisimple part is equal to $p$ reduces to the classification of nilpotent elements of the space $\tilde \z_1(p)$  with respect to the group $\tilde Z_0(p)$.
 
\subsection{Cartan subspace and the Weyl group} 
Any maximal subspace of $\g_1$ consisting of semisimple elements commuting with each other is called a Cartan subspace. All such spaces are transformed into each other by transformations from the group $G_0$.

Let $\c$ be some Cartan subspace. Let us put 
\[
N_0 (\c) = \{g = G_o\mid \forall x\in \c \;\; \g(x)\in \c\}, 
\]
\[
Z_0 (\c) = \{g = G_o\mid \forall x\in \c \;\; \g(x) = x\}, 
\]
\[
W = N_0(\c)/Z_0(\c).
\]
The group $W$ is finite. Considered as a group of linear transformations of the space $c$, it is called the Weyl group of the  graded Lie algebra $\g$ (with respect to $\c$). 

Every semisimple element $p \in \g_1$ can be translated by  a transformation from $G_0$ into a fixed Cartan subspace $\c$. Two elements from $\c$ are equivalent with respect to $G_0$ if and only if they are equivalent with respect to $W$. The restricted homomorphism
$\rho\colon \CC[\g_1]^{G_0} \to \CC[\c]^W$
of the algebra of invariants of the action of the group $G_0$ in the space $\g_1$ into the algebra of invariants of the action of the group $W$ in the space $\c$ is an isomorphism. 
The group $W$ is generated by complex reflections, which implies that the invariant algebra $\CC[\c]^W$ is free (see \cite{ShephardTodd54} or \cite{Chevalley55}).

\subsection{The support method}
Let $e$ be a non-zero nilpotent element of the space $\g_1$. 
Then, by Morozov's theorem, there exists a semisimple element $h\in \g_0$ and a nilpotent element $f\in \g_{-1}$ such that $[h, e] =2e$, $[h, f]=-2f, [e, f]= h$. Such an element $h$, which is defined up to conjugation by means of an element of the subgroup $Z_0(e)$, is called the characteristic of the element $e$. Let $\n_0$ denote the normalizer of the element $e$ in the algebra
 $\g_0$ (i.e. $\n_0 = \{ x \in \d_0 \mid [x,e] = ce\}$), with $\h$ its maximum torus. 
Denote by $\eta$ the linear form on $\h$ defined by the condition <<$ [u, e] = \eta(u).e$ for $u \in \h$>>, and consider the following $\ZZ$-graded subalgebra of the algebra $\g$:
\[
\overline\g  = \bigoplus \overline\g_k, \quad \overline\g_k = \{ x\in \g_{k\; \mathrm{mod}\, m} \mid [u,x] = k\eta(u) x \; \forall u \in \h \} 
\]

The algebra $\overline \g$ is reductive; its commutant considered as a $\ZZ$-graded subalgebra of $\g$ is called the support of the element $e$. The nilpotent element $e\in \g_1$ is determined by its support uniquely up to equivalence with respect to $G_0$. The support method is described in more detail in \cite{Vinberg75}.

\section{Semisimple four-vectors}\label{sec:3}
\subsection{Cartan subspace and the Weyl group}\label{sec:3.1}
Let us apply the results formulated above to our case, taking $\g$ as the algebra  of $E_7$ graded modulo 2 that was constructed in \S~\ref{sec:e7}, and take $G$ a simply connected Lie group having $\g$ as its Lie algebra. As already noted, in this case the orbits of the group $G_0$ in the space $\g_1 = \bw{4}V$ coincide with the orbits of the natural representation in this space of the group $\SL(V)$, which we and want to classify. 

In this situation, we indicate a specific Cartan subspace. 

An arbitrary permutation $s= (i_1 j_1 k_1 l_1 i_2 j_2 k_2 l_2)$ of the numbers $1,\ldots, 8$ can be associated with a semisimple four-vector 
\begin{equation}\label{eq:ps}
p(s) = e_{i_1} \wedge e_{j_1} \wedge e_{k_1} \wedge e_{l_1} + e_{i_2} \wedge e_{j_2} \wedge e_{k_2} \wedge e_{l_2} .
\end{equation}
It is easy to check that if the permutations $s= (i_1 j_1 k_1 l_1 i_2 j_2 k_2 l_2)$ and$s= (i_1' j_1' k_1' l_1' i_2' j_2' k_2' l_2')$ are such that the sets $\{i_1, j_1, k_1, l_1\}$ and $\{i_1', j_1', k_1', l_1'\}$ have exactly 2 common elements (in this case we say that $s$ and $s'$ are transversal), then the semisimple four-vectors $p(s)$ and $p(s')$ corresponding to them commute with each other. 
Consider the following seven permutations:
\[\begin{matrix}
s_1 = (1234\; 5678), && s_2 = (1357\; 6824), && s_3 = (1562\; 8437), \\
s_4 = (1683\; 4752), && s_5 = (1845\; 7263), && s_6 = (1476\; 2385), \\
&& s_7 = (1728\; 3546).
\end{matrix}
\]
 They are pairwise transversal, and therefore the space $\c\subset \g_1$ spanned by four-vectors $p_1 = p(s_1),\ldots, p_7=p(s_7)$ consists of mutually commuting semisimple elements. And since the dimension of the subspace $\c$ is equal to the rank of the algebra $\g$, then $\c$ is a Cartan subalgebra in $\g$, and the desired Cartan subspace in $\g_1$.
 
 Denote by $W$ the Weyl group of the graded algebra with respect to the Cartan subspace $\c$. Since $\c$ is a Cartan subalgebra, then $W$ coincides with the <<ordinary>> Weyl group of a semisimple Lie algebra $g$ with respect to $\c$ (see \cite{Vinberg76}). The group $W$ contains 63 reflections: its order is 2903040, and the degrees of the basic invariants are $2, 6, 8, 10, 12, 14, 18$. The four-vectors 
$p_1,\ldots, p_7$ form an orthogonal basis with respect to the Cartan scalar multiplication basis of the Cartan subspace $\c$, and their scalar squares are equal to each other.

The Weyl group $W$ contains reflections with respect to $p_1,\ldots, p_7$, generating in it a Weyl subgroup of type $7A_1$.

The root system of the algebra $\g$ with respect to $\c$, as a subset of $\c$, has the form
\begin{equation}
\Sigma = \{\gamma_i - \gamma_j, \ga_i + \ga_j + \ga_k + \ga_l \}\quad (i,j,k,l \text{ pairwise distinct}),
\end{equation}
$\begin{array}{lrcl}\text{where} 
&\ga_1 &=& c(p_1+ p_2 + p_3 + p_4 + p_5 + p_6 + p_7), \\
&\ga_2 &=& c(p_1- p_2 + p_3 - p_4 - p_5 - p_6 + p_7) ,\\
&\ga_3 &=& c(p_1+ p_2 - p_3 + p_4 - p_5 - p_6 - p_7) ,\\
&\ga_4 &=& c(p_1- p_2 - p_3 - p_4 + p_5 + p_6 - p_7) ,\\
&\ga_5 &=& c(-p_1+ p_2 + p_3 - p_4 + p_5 - p_6 - p_7), \\
&\ga_6 &=& c(-p_1- p_2 + p_3 + p_4 - p_5 + p_6 - p_7), \\
&\ga_7 &=& c(-p_1+ p_2 - p_3 - p_4 - p_5 + p_6 + p_7), \\
&\ga_8 &=& c(-p_1- p_2 - p_3 + p_4 + p_5 - p_6 + p_7) ,
\end{array}
$

\noindent with $c$ some arbitrary constant.

\subsection{Classification of semisimple four-vectors} It follows from the results formulated in \S~\ref{sec:2} that the classification of semisimple four-vectors reduces to the classification of elements of the Cartan subspace $\c$ constructed in \S~\ref{sec:3.1} with respect to the Weyl group $W$. 

Let $\Sigma$ be the root system of $\g$ with respect to the Cartan subalgebra $\c$. The centralizer $\z(p)$ of an arbitrary four-vector $p \in \c$ in the algebra $\g$ spans $\c$ and the root vectors corresponding to the roots that vanish on $p$; let $W_p$ be the Weyl subgroup of $W$ generated by the reflections corresponding to the roots that vanish on $p$. It is known that $W_p$ coincides with the stabilizer subgroup of the element $p$ in the group $W$ (this is one of the properties of groups generated by reflections). Thus, the stabilizer subgroup $W_p$, and the centralizer $\z(p)$ define each other uniquely. Let's put 
\[\begin{matrix}
\c_p=\{x\in \c \mid \forall w \in W_p\;\; w (x) = x\}, \\
\c^0_p = \{x \in \c: W_x = W_p\}. 
\end{matrix}
\]
If for $p,q \in \c$, $\c_p = w (\c_q)$, then $W_q = wW_pw^{-1}$, and vice versa. In particular, the subspace $\c_p$ is invariant under the normalizer $N(W_p)$ of the subgroup $W_p$ in the group $W$. The quotient group $\Gamma_p = N(W_p)/W_p$ acts in a natural way on the space $\c_p$; two elements of the subset $\c^0_p \subset \c_p$ are equivalent with respect to $W$ if and only if they are equivalent with respect to $\Gamma_p$. 

The group $W$ has 32 classes of conjugate subgroups $W_p$.  
Table~\ref{table:1}  shows for each class the orders of the groups $W_p$ and $\Gamma_p$ (calculated using the tables in \cite{ShephardTodd54}), the type of the subalgebra $\tilde \z (p) = \z (p)' $ is indicated, and also for the representative of this class the basis of the subspace $\c_p$. The following abbreviations are used: the four-vector $p_{k_1} + \cdots + p_{k_n}$ is denoted by the symbol $P(k_1\ldots k_n)$. 

A semisimple four-vector $p$ will be called canonical if $p\in \c$ and its subspace $\c_p$ coincides with one of the 32 subspaces whose bases are given in Table~\ref{table:1} . The semisimple four-vectors can be divided into 32 families, referring a semisimple four-vector to the $k$-th family if the stabilizer subgroup (in $W$) of its equivalent canonical four-vector (it does not depend on the choice of the canonical representative of the orbit of the given four-vector) belongs to $k$-th class of conjugated subgroups $W_p$. Arbitrary (non-semisimple) four-vectors will be assigned to those families to which their semisimple parts belong.

\begin{table}
\caption{The classes of semi-simple elements in $\bw{4}\CC^8$.}\label{table:1}
\begin{tabular}{r | l | c | r | r }
\hline
\textnumero & Basis $\c_p$ & Type $\z_p'(a)$  & $\left| W_p \right|$ & $\left| \Gamma_p \right|$\\  \hhline{=|=|=|=|=}
 1 && $E_7$ & 2\;903\;040 & 1 \\  \hline
 2 & $P(1)+P(3)-P(7)$ & $E_6$ & 51\;840 & 2 \\  \hline
 3 & $P(1) $& $D_6$ & 23\;040 & 2 \\  \hline
 4 & $P(3) -P(7)$& $D_5 + A_1$ & 3\;840 & 2 \\  \hline
 5 & $P(1\;234\;567)$ & $A_6$ & 5\;040 & 2 \\  \hline
 6 & $P(123\;456) +P(45)$ & $A_5+A_1$ & 1\;440 & 2 \\  \hline
 7 & $P(123\;456) +P(146)$ & $A_4+A_2$ & 720 & 2 \\  \hline
 8 & $P(123\;456) $ & $A_3+A_2+A_1$ & 288 & 2 \\  \hline
 9 & $P(1), P(13)-P(7)$ & $D_5$ & 1\;920 & 4 \\  \hline
 10 & $P(1), P(3) $& $D_4 + A_1$ & 384 & 8 \\  \hline
 11 & $P(2), P(13)-P(7)$ & $A_5$ & 720 & 4 \\  \hline
 12 & $P(1), P(26)-P(3) $& $A_5$ & 720 & 12 \\  \hline
 13 & $P(2\;345)+P(4), P(156) $& $A_4+A_1$ & 240 & 2 \\  \hline
 14 & $P(1\;235)+P(1), P(1\;346)+P(1) $& $A_3+A_2$ & 144 & 4 \\  \hline
 15 & $P(13), P(2\;456) $& $A_3+2A_1$ & 96 & 4 \\  \hline
 16 & $P(156), P(23)+2P(4)+ 3P(5) $& $2A_2+A_1$ & 72 & 4 \\  \hline
 17 & $P(123\;456), P(1)-P(3) $& $A_2+3A_1$ & 48 & 12 \\  \hline
 18 & $P(1), P(3),P(7) $& $D_4$ & 192 & 48 \\  \hline
 19 & $P(1), P(13)-P(7), P(26)-P(3)$ & $A_4$ & 120 & 12 \\  \hline
 20 & $P(13), P(2\;456), P(7)$ & $A_3+A_1$ & 48& 8 \\  \hline
 21 & $P(13), P(26), P(45)$ & $A_3+A_1$ & 48& 48 \\  \hline
 22 & $P(1\;236), P(15)-P(6), P(124)+2P(5) $& $2A_2$ & 36& 24 \\  \hline
 23 & $P(5), P(16), P(23)+2P(4)$ & $A_2+2A_1$ & 24& 8 \\  \hline
 24 & $P(1), P(2), P(3) $& $4A_1$ & 16& 48 \\  \hline
 25 & $P(1), P(56), P(1\;236),P(24)$ & $A_3$ & 24& 96 \\  \hline
 26 & $P(5), P(16), P(246),P(143)-P(5)$ & $A_2+A_1$ & 12& 48 \\  \hline
 27 & $P(1), P(3), P(3),P(6) $& $3A_1$ & 8& 1\;152 \\  \hline
 28 & $P(1), P(2), P(3),P(4) $& $3A_1$ & 8& 96 \\  \hline
 29 & $P(1), P(5), P(6),P(24),P(23) $& $A_2$ & 6& 1\;140 \\  \hline
 30 & $P(1), P(2), P(3),P(4),P(5)$& $2A_1$ & 4& 768 \\  \hline
 31 & $P(1), P(2), P(3),P(4),P(5),P(6)$& $A_1$ & 2& 23\;040 \\  \hline
 32 & $P(1), P(2), P(3),P(4),P(5),P(6),P(7)$&  & 1& 2\;903\;040 \\  \hline
\end{tabular}
\end{table}

\section{Classification of Nilpotent Four-Vectors}\label{sec:4}
Table~\ref{table:2} shows the results of the classification of nilpotent four-vectors obtained by the support method applied to the graded Lie algebra $\g$  constructed in \S~\ref{sec:an12}. 

In each equivalence class of nilpotent four-vectors in the corresponding row of the table, the <<type>> indicates the support type of the elements of this class is indicated as $\ZZ$-graded subalgebras of the graded algebra $\g$ (in the notation of \cite{Vinberg75}). For each class a representative is chosen that has a characteristic that is contained in the Cartan subalgebra $\d$ of $\g_0$ and whose numerical marks are relative to the system of simple roots $\Pi_0 = \{ \ep_1- \ep_2, \ep_2-\ep_3,\ldots, \ep_7-\ep_8\}$ are non-negative; they are indicated in the <<Characteristics>> column. In addition, the dimension $d$ of the stabilizer subalgebra of the chosen representative and the type of its maximal reductive subalgebra $\mathfrak{s}_0$ are given.

If in Table~\ref{table:2}, the number $i$ is the class of the nilpotent element $e\in \g_1$, which has characteristic $h \in \g_0$, and the number $i+1$ is omitted, this means that $i+1$ is the number of the class of the nilpotent element $f \in \g_{-1}= \g_1$, which together with $h$ and $e$ forms a standard basis for a simple three-dimensional subalgebra and hence has characteristic $-h$.

\begin{table}
\caption{The classes of nilpotent elements in $\bw{4}\CC^8$. \emph{LO: We have corrected column $d$ striking the original and replacing it in italics with the correct value.  }
}\label{table:2}
\resizebox{1 \textwidth}{!}{\begin{tabular}{r | c | c | r | c  ||r | c | c | r | c  }
\hline
\textnumero  & Characteristic & Type  & $d$ & Type $\s_0$ &\textnumero  & Characteristic & Type  & $d$ & Type $\s_0$ \\  
\hhline{=|=|=|=|=||=|=|=|=|=}
$1$ & $\ch0001000$ & $A_1$ & $46$ & $2A_3$          &		$44$ & $\ch0020200$ & $D_4(a_1)$ & $16$ &$ T_3 $ \\
$2$ & $\ch0100010$ & $2A_1$ & $37$ & $C_2+ 2A_2+T_1$  &	 $45$ & $\ch2020020$ & $A_3+ A_2+ A_1$ & $13$ &$ T_1 $ \\
$3$ & $\ch0200000$ & $3A_1$ & $36$ & $C_3+ A_1$          &    $47$ & $\ch0040000$ & $A_3+ A_2+ A_1$ & $13$ &$ A_1 $ \\
$5$ & $\ch1001001$ & $3A_1$ & $31$ & $A_2+ T_2$          &    $49$ & $\ch1111111$ & $2A_3+ A_1$ & \st{$10$}\; \emph{11} &$ 0 $ \\
$6$ & $\ch0002000$ & $4A_1$ & $30$ & $A_3$          &       	$50$ & $\ch2004002$ & $D_4$ & $15$ &$ A_2+T_1 $ \\
$7$ & $\ch1100100$ & $4A_1$ & $28$ & $A_2+ T_1$          &    $51$ & $\ch2204022$ & $D_5$ & $7$ &$ T_2 $ \\
$9$ & $\ch2000002$ & $A_2$ & $30$ & $A_2+ A_2 + T_1$  &    $52$ & $\ch1313143$ & $D_6$ & $4$ &$ T_1 $ \\
$10$ & $\ch2010001$ & $A_2+A_1$ & $25$ & $2A_1+T_2$  &    $54$ & $\ch0202202$ & $D_5(a_1)+A_1$ & $9$ &$ T_1 $ \\
$12$ & $\ch0101010$ & $5A_1$ & $25$ & $2A_1$  &    		$56$ & $\ch0004004$ & $D_5(a_1)+A_1$ & $9$ &$ T_1 $ \\
$13$ & $\ch3000100$ & $A_2+ 2A_1$ & $22$ & $2A_1+T_1$  & $58$ & $\ch1013012$ & $D_4+A_1$ & $12$ &$ A_1 + T_1 $ \\
$15$ & $\ch1010101$ & $A_2+ 2A_1$ & $22$ & $T_3$  & $60$ & $\ch1112111$ & $D_4+2A_1$ & $10$ &$ T_1 $ \\
$16$ & $\ch4000000$ & $A_2+ 3A_1$ & $21$ & $G_2$  & $61$ & $\ch0103103$ & $D_5(a_1)$ & $10$ &$ T_2 $ \\
$18$ & $\ch2000200$ & $A_2+ 3A_1$ & $21$ & $2A_1$  & $63$ & $\ch3113121$ & $D_5 + A_1$ & $6$ &$ T_1 $ \\
$20$ & $\ch0102010$ & $A_3$ & $21$ & $3A_1+T_1$  & $65$ & $\ch0404044$ & $E_7(b)$ & $3$ &$ 0 $ \\
$21$ & $\ch2200022$ & $A_4$ & $13$ & $A_1+T_2$  & $67$ & $\ch2422222$ & $D_6+A_1$ & $3$ &$ 0 $ \\
$22$ & $\ch1211121$ & $A_5$ & \st{$7$}\;\emph{9} & $T_2$  & $69$ & $\ch4224224$ & $E_6$ & $3$ &$ T_1 $ \\
$23$ & $\ch0202040$ & $A_5$ & $12$ & $2A_1$  & $70$ & $\ch3013131$ & $D_6(a_1)$ & $6$ &$ T_1 $ \\
$25$ & $\ch2220222$ & $A_6$ & $6$ & $T_1$  & $72$ & $\ch1110111$ & $D_4(a_1)+2A_1$ & \st{$11$}\;\emph{14} &$ T_1 $ \\
$26$ & $\ch2222222$ & $A_7$ & $4$ & $0$  & $73$ & $\ch1010210$ & $D_4(a_1)+A_1$ & $15$ &$ T_2 $ \\
$27$ & $\ch0200020$ & $A_2+A_2$ & $21$ & $2A_1+T_1$  & $75$ & $\ch2022222$ & $D_6(a_1)+A_1$ & $5$ &$ 0 $ \\
$28$ & $\ch1030010$ & $A_2+A_3$ & $14$ & $T_2$  & $77$ & $\ch4004040$ & $E_7(c_1)$ & $5$ &$ 0 $ \\
$30$ & $\ch0040040$ & $E_7(c_2)$ & $7$ & $0$  & $79$ & $\ch4220224$ & $E_6(a_1)$ & $4$ &$ T_1 $ \\
$32$ & $\ch2020202$ & $A_2+A_4$ & $10$ & $T_1$  & $80$ & $\ch1030131$ & $D_6(a_2)$ & $8$ &$ T_1 $ \\
$33$ & $\ch0202020$ & $A_3+A_3$ & $13$ & $A_1$  & $82$ & $\ch2020222$ & $D_6(a_2)+A_1$ & $7$ &$ 0 $ \\
$34$ & $\ch0002020$ & $A_3+A_1$ & $20$ & $3A_1$  & $83$ & $\ch4444044$ & $E_7(a_1)$ & $1$ &$ 0 $ \\
$36$ & $\ch1011101$ & $A_3+A_1$ & \st{$14$}\;\emph{17} & $T_3$  & $86$ & $\ch4404404$ & $E_7(a_2)$ & $2$ &$ 0 $ \\
$37$ & $\ch3101021$ & $A_4+A_1$ & $11$ & $T_2$  & $88$ & $\ch4444448$ & $E_7$ & $0$ &$ 0 $ \\
$39$ & $\ch2202022$ & $A_5+A_1$ & $8$ & $T_1$  & $90$ & $\ch0101111$ & $A_3+2A_1$ & $16$ &$ T_2 $ \\
$40$ & $\ch1311111$ & $A_5+A_1$ & \st{$8$}\;\emph{9} & $T_1$  & $92$ & $\ch2002002$ & $A_3+2A_1$ & $16$ &$ A_1  + T_1 $ \\
$42$ & $\ch0220220$ & $E_6(b)$ & $8$ & $T_1$  & $93$ & $\ch2101101$ & $A_3+3A_1$ & $15$ &$ A_1  $ \\
$43$ & $\ch1101011$ & $2A_2+A_1$ & \st{$15$}\;\emph{18} & $T_2$  & $94$ & $\ch1011012$ & $A_3+3A_1$ & $15$ &$ A_1   $\\
\end{tabular}
}
\end{table}

\section{Mixed Four-Vectors}\label{sec:5}

\subsection{Some preliminary remarks} We call a four-vector mixed if it is neither semisimple nor nilpotent. The semisimple part of an arbitrary mixed four-vector belongs to one of the 30 families appearing in Table~\ref{table:1} numbered 2-31. It follows from the results formulated in \S~\ref{sec:2} that the classification of four-vectors of a fixed family reduces to the classification of four-vectors with a semisimple part lying in the corresponding subspace $\c_p$, whose basis is given in Table~\ref{table:1} (see \S~\ref{sec:3}). The classification of four-vectors with a given semisimple part $p$ reduces, in turn (in the notation of \S~\ref{sec:2.1}), to the classification of nilpotent elements in the subspace $\tilde \z_1(p)$ of the graded Lie algebra $\tilde \z(p)$ with respect to the group $\tilde{Z}_0(p)$. 
By the support method one can classify the nilpotent four-vectors of the algebra $\tilde{\z}(p)$ with respect to the connected component of the identity of the group $\tilde{Z}_0(p)$. 
Then we need to check which of the resulting classes are equivalent with respect to groups of connected components. 

In the following sections we give a classification of the nilpotent parts of mixed four-vectors of some of the most interesting families. 

In each case, in the algebra $\tilde\z_0(p)$ one chooses a Cartan subalgebra and a system of simple roots with respect to this Cartan subalgebra is indicated. The characteristics of nilpotent parts
of the four-vectors of the corresponding family are taken to be contained in the chosen Cartan subalgebra and dominant; their corresponding numerical marks are indicated in the  <<Characteristics>> column of the corresponding table. The remark made at the end of \S~\ref{sec:4} still holds. 

The notation introduced in \S~\ref{sec:e7} is used.

\subsection{Mixed four-vectors, second family}\label{sec:5.2}
The canonical four-vector $p=2(p_1+p_3-p_7)$ belongs to the second family. In the exterior algebra of the space $V$ it is equal to the square of the non-degenerate bivector $e_1\wedge e_2 +e_3\wedge e_4 +e_5\wedge e_6 +e_7\wedge e_8 $, so the connected component of the identity of its stabilizer subgroup $S(p)$ in the group $\SL(V)$ is a simple group of type $C_4$. It consists of all unimodular transformations of the space $V$ that preserve the non-degenerate skew-symmetric bilinear functional 
\[
B=e^1\wedge e^2 + e^3\wedge e^4 + e^5 \wedge e^6  + e^7\wedge e^8.
\]
The semisimple component $\g^{(2)}$ of the centralizer of the four-vector $p$ in algebra $\g$ is a graded simple subalgebra of type $E_6$ (see Table~\ref{table:1}), and
\begin{equation}
\g^{(2)}_0 = \{ X \in L_0(V) \mid X^s_i B_{sj} + X^s_j B_{is} = 0\},
\end{equation}
\begin{equation}
\g^{(2)}_1 = \{ T \in \bw4V \mid  T^{srij} B_{sr} = 0\}.
\end{equation}

Table~\ref{table:3} lists the classes of nilpotent four-vectors of the subalgebra $\g^{(2)}$. The Cartan subalgebra of $\g_0^{(2)}$ is taken as the subspace $\d^{(2)} = \d \cap \g_0^{(2)}$. It is easy to check that 
\[
\d^{(2)} = \{A \in \d\mid (\ep_1 + \ep_2)(A) = (\ep_3 + \ep_4)(A) = (\ep_5 + \ep_6)(A) = (\ep_7 + \ep_8)(A) = 0 \}.
\]
The following are chosen as simple roots of the algebra $\g_0^{(2)}$ with respect to $\d^{(2)}$: 
\[
\ep_1 - \ep_3,\; \ep_3-\ep_5,\; \ep_5-\ep_7,\; \ep_7-\ep_8.
\]
\begin{table}
\caption{Nilpotent parts of four-vectors of the 2nd family}\label{table:3}
\begin{tabular}{r | c | c || r | c | c   }
\hline
\textnumero  & Characteristic & Type  &\textnumero  & Characteristic & Type   \\  \hhline{=|=|=||=|=|=}
$1$ & $\ch0001$ & $A_1$ & 			$13$ & $\ch1111$ &  $A_3 + A_1$         	 \\
$2$ & $\ch0100$ & $A_1$ & 			$14$ & $\ch2004$ &  $B_3 $         	 \\
$3$ & $\ch1001$ & $2A_1$ & 			$15$ & $\ch1211$ &  $C_3$         	 \\
$4$ & $\ch0002$ & $2A_1$ & 			$16$ & $\ch2222$ &  $C_4 $         	 \\
$5$ & $\ch2000$ & $A_2$ & 			$17$ & $\ch2204$ &  $F_4(e_1) $         	 \\
$6$ & $\ch0101$ & $3A_1$ & 			$18$ & $\ch1011$ &  $C(1,2)$         	 \\
$7$ & $\ch1010$ & $A_2+A_1$ & 		$19$ & $\ch2002$ &  $C_2+A_1$         	 \\
$8$ & $\ch0200$ & $A_2$ & 			$20$ & $\ch0220$ &  $F_4(e_2)$         	 \\
$9$ & $\ch0102$ & $C_2$ & 			$21$ & $\ch2202$ &  $C_3 + A_1$         	 \\
$10$ & $\ch0202$ & $A_3$ & 			$22$ & $\ch1112$ &  $B_3 + A_1$         	 \\
$11$ & $\ch0020$ & $F_4(e_3)$ & 		$23$ & $\ch4224$ &  $F_4$         	 \\
$12$ & $\ch1101$ & $A_2+A_1$  			
\end{tabular}
\end{table}

\subsection{Mixed four-vectors, third family}\label{sec:5.3} Let us suppose that a decomposition of the space $V$ into a direct sum of two four dimensional subspaces is given: 
\begin{equation}\label{eq5.3}
V  = W_1 \oplus W_2.
\end{equation}
Then we can identify the space $(\bw{2} W_1) \otimes (\bw2W_2)$ with the subspace in $\bw4V$, assuming that 
\[
(w_1\wedge w_1' )\otimes (w_2 \wedge w_2') = w_1\wedge w_1' \wedge w_2 \wedge w_2' \quad (w_k, w_k' \in W_k)
.\]
The algebra $L_0(W_1) \oplus L_0(W_2)$ is naturally identified as a subalgebra of the algebra $L_0(V)$;
it is easy to check that the subspace
\[\g^{(3)} = (L_0(W_1) \oplus L_0(W_2)) \oplus ((\bw2W_1 )\otimes( \bw2W_2))
\]
 is a graded subalgebra of the graded algebra $\g$.
 
 Suppose that the vectors $e_1, e_2, e_3, e_4$  form the basis of the subspace $W_1$, and the vectors $e_5, e_6, e_7, e_8$  form the basis of the subspace
 $W_2$, i.e. $W_1 = \langle e_1, e_2, e_3, e_4 \rangle$, $W_2 = \langle e_5, e_6, e_7, e_8 \rangle$). 
 Then $\g^{(3)}$ is a graded subalgebra with respect to $\d$ of the algebra $\g$. An investigation of the system of roots shows that $\g^{(3)}$ is a simple algebra of type $D_6$, and $\g^{(3)}_0 = \g^{(3)}  \cap \g_0$  is its semisimple subalgebra of type $2A_3$. 
 
It is easy to check that $\g^{(3)}$ is contained in the centralizer of the semisimple four-vector $p=p_1$. This four-vector belongs to the third family, so that the semisimple component of its centralizer is a simple subalgebra of type $D_6$ (see Table~\ref{table:1}) and hence coincides with $\g^{(3)}$.

The stabilizer subgroup $S(p)$ of the four-vectors $p$ in the group $\SL(V)$ is an extension of the group $\SL(W_1)\times \SL(W_2)$ by a group of order $2$ that permutes the subspaces $W_1$ and $W_2$. Table~\ref{table:4} lists the classes of nilpotent four-vectors of the subalgebra $\g^{(3)}$. The subspace $\g^{(3)}_0$  is taken as the Cartan subalgebra $\d^{(3)} = \g^{(3)}_0\cap \d$ (it is simultaneously a Cartan subalgebra of the algebra $\g^{(3)}$). 
Obviously, 
\[
\d^{(3)} = \{ A\in \d \mid (\ep_1 +\ep_2 + \ep_3 + \ep_4)(A) = 0\}.
\]
The following are taken as simple roots of the algebra $\g^{(3)}_0$ with respect to $\d^{(3)}$:
\[
\ep_1 - \ep_2,\quad \ep_2-\ep_3, \quad\ep_3-\ep_4,\quad \ep_5-\ep_6,\quad  \ep_6-\ep_7,\quad  \ep_7-\ep_8.  
\]

\begin{table}
\caption{Nilpotent parts of four-vectors of the 3rd family}\label{table:4}
\begin{tabular}{r | c | c || r | c | c   }
\hline
\textnumero  & Characteristic & Type  &\textnumero  & Characteristic & Type   \\  \hhline{=|=|=|| |=|=|=}
$1$ & $\ch010010$ & $A_1$ & 				$25$ & $\ch224224$ &  $A_5$         	 \\
$2$ & $\ch020000$ & $2A_1$ & 			$27$ & $\ch222222$ &  $2A_3$         	 \\
$3$ & $\ch101101$ & $2A_1$ & 			$28$ & $\ch040242$ &  $D_4$         	 \\
$4$ & $\ch111010$ & $3A_1$ & 			$29$ & $\ch444242$ &  $D_5$         	 \\
$5$ & $\ch200200$ & $3A_1$ & 			$30$ & $\ch444844$ &  $D_6$         	 \\
$7$ & $\ch200002$ & $3A_1$ & 			$32$ & $\ch131111$ &  $A_3+A_1$         	 \\ 
$8$ & $\ch210101$ & $4A_1$ & 			$33$ & $\ch220022$ &  $A_3+A_1$         	 \\
$10$ & $\ch202000$ & $A_2$ & 			$34$ & $\ch220220$ &  $A_3+A_1$         	 \\
$11$ & $\ch020020$ & $4A_1$ & 			$36$ & $\ch131341$ &  $D_4+A_1$         	 \\
$12$ & $\ch301010$ & $A_2+A_1$ & 		$38$ & $\ch222242$ &  $D_4+2A_1$         	 \\
$14$ & $\ch111111$ & $5A_1$ & 			$39$ & $\ch040440$ &  $D_5(a_1)$         	 \\
$15$ & $\ch400000$ & $A_2+2A_1$ & 		$41$ & $\ch222020$ &  $A_3+2A_1$         	 \\
$17$ & $\ch202020$ & $A_2+2A_1$ & 		$42$ & $\ch040202$ &  $D_4(a_1)$         	 \\
$18$ & $\ch202202$ & $A_2+A_2$ & 		$43$ & $\ch311121$ &  $A_3+2A_1$         	 \\
$19$ & $\ch040400$ & $A_3+A_2$ & 		$45$ & $\ch044444$ &  $D_6(a_1)$         	 \\
$21$ & $\ch040020$ & $A_3$ & 			$47$ & $\ch044404$ &  $D_6(a_2)$         	 \\
$22$ & $\ch121121$ & $A_3$ & 			$49$ & $\ch103131$ &  $D_4(a_1) + A_1$         	 \\
$23$ & $\ch404202$ & $A_4$ & 			$51$ & $\ch202222$ &  $D_4(a_1)+ 2A_1$         	 \\
$24$ & $\ch224422$ & $A_5$ & 			
\end{tabular}
\end{table}

\subsection{Mixed Four-Vectors, Ninth family}\label{sec:5.4} Let  $p$ be the canonical semisimple four-vector of the ninth family, i.e., $p\in \c$, and the subspace $\c_p$ coincides with the linear span of the four-vectors $p_1$ and $p_1+p_3-p_7$. These latter are the canonical semisimple four-vectors of the third and second families, respectively, so the semisimple component $\g^{(9)}$ of the centralizer of the torus of the four-vector $p$ in $\g$ (which is a semisimple algebra of type $D_5$ - see Table~\ref{table:1}) coincides with the intersection of subalgebras $\g^{(2)}$ and $\g^{(3)}$ constructed in Sections~\ref{sec:5.2} and \ref{sec:5.3}, i.e. 
\begin{equation}
\g_0^{(9)} = \{ X \in L_0(W_1) \oplus L_0(W_2) \mid X^s_i B_{sj} + X^s_jBV_{is} = 0\} \subset \g_0,
\end{equation}
\begin{equation}
\g_1^{(9)} = \{ T \in (\bw2 W_1) \otimes (\bw2 W_2) \mid T^{srij}B_{sr} = 0\} \subset \g_1,
\end{equation}
with $W_1 = \langle e_1, e_2, e_3, e_4\rangle$, $W_2 = \langle e_5, e_6, e_7, e_8 \rangle$,
\[
B = e^1 \wedge e^2 + e^3\wedge e^4 + e^5 \wedge e^6 + e^7\wedge e^8 \in \bw2 V^*.
\]
In the algebra $\g^{(9)}$ the subalgebra $\g_0^{(9)}$ is a semisimple subalgebra of type $2C_2$: it is naturally isomorphic to the direct sum of the stabilizer subalgebra in $L_0(W_1)$ of the non-degenerate bivector  $e^1\wedge e^2 + e^3\wedge e^4\in \bw2 W_1^*$ and the stabilizer subalgebra in $L_0(W_2)$ of the non-degenerate bivector $e^5\wedge e^6 + e^7 \wedge e^8 \in \bw2W_2^*$.

The group of connected components of the stabilizer subgroup $S(p)$ is not trivial; it is essential for us that there are transformations in $S(p)$ that permute the subspaces $W_1$ and $W_2$: such transformations induce outer automorphisms of the stabilizer subalgebras  $\z(p)$ of four-vectors $p$ in algebra $\g_0$.

Table~\ref{table:5} lists the classes of nilpotent four-vectors of subalgebras of $\g^{(9)}$. The subspace $\d^{(9)} = \d \cap \g^{(9)}$ is taken as a Cartan subalgebra in $\g^{(9)}_0$. Obviously, $\d^{(9)} = \d^{(2)}$. 
The following are taken as simple roots of the algebra $\g^{(9)}_0$ with respect to $\d^{(9)}$:
\[
\ep_1 - \ep_3,\quad \ep_3-\ep_4, \quad\ep_5-\ep_7,\quad \ep_7-\ep_8.
\]

\begin{table}
\caption{Nilpotent parts of four-vectors of the 9th family}\label{table:5}
\begin{tabular}{r | c | c || r | c | c   }
\hline
\textnumero  & Characteristic & Type  &\textnumero  & Characteristic & Type   \\  \hhline{=|=|=||=|=|=}
$1$ & $\ch0101$ & $A_1$ & 			$10$ & $\ch1212$ &  $A_3$         	 \\
$2$ & $\ch0200$ & $A_1$ & 			$11$ & $\ch0420$ &  $B(0,1,2)$         	 \\
$3$ & $\ch1010$ & $2A_1$ & 			$12$ & $\ch2404$ &  $B_3$         	 \\
$4$ & $\ch1101$ & $2A_1$ & 			$13$ & $\ch2444$ &  $B_4$         	 \\
$5$ & $\ch0202$ & $2A_1$ & 			$14$ & $\ch1113$ &  $C_2+A_1$         	 \\
$6$ & $\ch2000$ & $A_2$ & 			$15$ & $\ch0222$ &  $C_2+A_1$         	 \\
$7$ & $\ch1111$ & $3A_1$ & 			$16$ & $\ch2222$ &  $C_2+C_2$         	 \\
$8$ & $\ch2002$ & $A_2+A_1$ & 			$17$ & $\ch2422$ &  $B_3+A_1$         	 \\
$9$ & $\ch0402$ & $C_2$ & 		
\end{tabular}
\end{table}

\subsection{Mixed Four-Vectors, Eleventh Family}\label{sec:5.5}
Let us assume that in addition to decomposition \eqref{eq5.3} the following  isomorphisms are given
\begin{equation}\label{eq5.7}
\phi_1\colon W \to W_1,\quad \phi_2\colon W \to W_2,
\end{equation}
where $W$ is some four-dimensional vector space.

Let $\{f_1,f_2, f_3, f_4\}$ be a basis of $W$. We assume that 
\[
\phi_1(f_k) = e_{2k-1},\quad \phi_2(f_k) = e_{2k} \quad (k=1,2,3,4).
\]

Let $\rho$ be an embedding of the algebra $L_0(W)$ into the algebra $L_0(V)$, which associates an element of the algebra $L_0(W)$ written in the basis $\{f_1,f_2, f_3, f_4\}$ as the matrix $A$ with a transformation of the space $V$ written in the basis 
$\{\phi_1(f_1),\ldots, \phi_1(f_4), \phi_2(f_1),\ldots, \phi_2(f_4)\}$ as the matrix $\begin{pmatrix} A & 0 \\ 0 & -A^*\end{pmatrix}$,  (<<$^*$>> means transposition).

The subalgebra $\g^{(11)}_0 =\rho(L_0(W))$ of the algebra $L_0(V)$ is contained in the subalgebra $L_0(W_1) \oplus L_0(W_2)$:
\[
\g^{(11)}_0 = \{ X\in L_0(W_0) \oplus L_0(W_2) \subset L_0(V) \mid X^s_i B_{sj} + X^s_j B_{is} = 0\},
\]
where $B = e^1\wedge e^2 + e^3\wedge e^4 + e^5\wedge e^6 + e^7 \wedge e^8$.

Considering $\bw2W_1 \otimes \bw2W_2$ as a subspace in $\bw4V$, we set 
\[
\g^{(11)}_0 = \{T \in \bw2W_1\otimes \bw2W_2 \mid T^{srij}B_{sz} = 0  \} .
\]
Let $p$ be the canonical semisimple four-vector of the 11th family, i.e., $p\in \c$ and $\c_p = \langle p_1+p_3- p_7, p_2\rangle$. Then the semisimple component of the centralizer of the four-vector $p$ in the algebra $\g$ (which is a subalgebra of type $A_5$) coincides with the subalgebra 
\[\g^{(11)} = \g^{(11)}_0 \oplus \g^{(11)}_1. 
\]
The stabilizer subgroup $S(p)$ of a four-vector in the group $SL(V)$ is not connected: with the help of its elements one can realize any automorphism of the subalgebra $\g^{(11)}$. 

Table~\ref{table:6} lists the classes of nilpotent four-vectors of the subalgebra $\g^{(11)}$. The subspace $\d^{(11)} = \d \cap \g_0^{(11)}$ is taken as the Cartan subalgebra in $\g^{(11)}_0$. It is easy to check that 
\[
\d^{(11)} = \{ A\in \d^{(2)} \mid (\ep_1 + \ep_3 + \ep_5 + \ep_7)(A)= 0\} .
\]

The following are taken as simple roots of the algebra $\g^{(11)}_0$ with respect to $\d^{(11)}$:
\[
\ep_1 - \ep_3,\quad \ep_3-\ep_5, \quad\ep_5-\ep_7.
\]

\begin{table}
\begin{tabular}{r | c | c || r | c | c   }
\hline
\textnumero  & Characteristic & Type  &\textnumero  & Characteristic & Type   \\  \hhline{=|=|= =|=|=}
$1$ &  {}$\ch010${} & $A_1$ & 			$6$ & {}$\ch202${}&  $A_2$         	 \\
$2$ &  {}$\ch101${} & $A_1$ & 			$7$ & {}$\ch222${} &  $A_3$         	 \\
$3$ &  {}$\ch200${} & $2A_1$ & 		$8$ & {}$\ch121${} &  $C_2$         	 \\
$4$ &  {}$\ch020${} & $2A_1$ & 		$9$ & {}$\ch022${} &  $C_2+A_1$         	 \\
$5$ &  {}$\ch111${} & $3A_1$ & 		$10$ & {}$\ch224${} &  $C_3$         	 
\end{tabular}
\caption{Nilpotent parts of four-vectors of the 11th family}\label{table:6}
\end{table}
\subsection{Mixed Four-Vectors, Twelfth Family}\label{sec:5.6}
Let decomposition \eqref{eq5.3} and isomorphisms \eqref{eq5.7} be given. Then the space $S^2(\bw2W)$ can be identified with a subspace in $\bw4V$, assuming that 
\begin{multline*}
(w_1\wedge w_1')(w_2\wedge w_2') = 
\phi_1(w_1)\wedge \phi_1(w_1') \wedge \phi_2(w_2)\wedge \phi_2(w_2')
+\\
\phi_2(w_1)\wedge \phi_2(w_1') \wedge \phi_1(w_2)\wedge \phi_1(w_2')
\quad (w_k, w_k' \in W).
\end{multline*}
With the isomorphisms $\phi_1$ and $\phi_2$ one can associate the embedding of the algebra $L_0(W)$ into the algebra $L_0(W_1)\oplus  L_0(W_2)$ (and thus into the algebra $L_0(V)$) in the form of a <<diagonal>> subalgebra. Let's put
\[\g^{(12)}_0 = L_0(W) \subset \g_0,
\]
\[\g^{(12)}_1 = \{ T\in S^2 (\bw2W) \mid T^{k_1l_1k_2l_2}d_{k_1l_1k_2l_2} = 0\} \subset \g_1,
\]
where $d$ is an arbitrary nonzero element of the space $\bw4W^*$, and  $S^2(\bw2W)$ is considered as a subspace of $\otimes^4W$. The subspace $\g^{(12)} = \g^{(12)}_0 \oplus \g^{(12)}_1$. 
is a graded simple subalgebra of type $A_5$ in $\g$. Let $\{f_1, f_2, f_3, f_4\}$ be a basis in $W$; we will assume that
\[
\phi_1(f_1) = e_1,\quad \phi_1(f_2) = e_2,\quad \phi_1(f_3) = e_3,\quad \phi_1(f_4) = e_4,
\]
\[
\phi_2(f_1) = e_8,\quad \phi_2(f_2) = e_7,\quad \phi_2(f_3) = e_6,\quad \phi_2(f_4) = e_5.
\]
Then an element of the algebra $L_0(W)$ written in the basis $\{f_1, f_2, f_3, f_4\}$ as the matrix $A$, in the basis $\{e_1, e_2, e_3, e_4, e_8, e_7, e_6, e_5\}$ is the matrix $\begin{pmatrix} A& 0 \\ 0 & A\end{pmatrix}$.

Let $p$ be the canonical semisimple four-vector of the 12th family, i.e., set $p \in \c$ and $\c_p  = \langle p_1, p_2-p_3+p_6\rangle$. Then the commutant and its centralizer in the algebra $\g$ coincide with the one constructed above under the algebra $\g^{(12)}$. The stabilizer subgroup $S(p)$ is not connected, but its elements induce inner automorphisms of the subalgebra $\g^{(12)}_0$.

Table~\ref{table:7} lists the classes of nilpotent four-vectors of the subalgebra $\g^{(12)}$. The subspace $\d^{(12)} = \d \cap \g^{(12)}_0$ is taken as the Cartan subalgebra in $\g^{(12)}_0$. 

The following are taken as simple roots of the algebra $\g^{(12)}_0$ with respect to $\d^{(12)}$:
\[
\ep_1 - \ep_2,\quad \ep_2-\ep_3, \quad\ep_3-\ep_4
.\]
\begin{table}
\caption{Nilpotent parts of four-vectors of the 12th family}\label{table:7}
\begin{tabular}{r | c | c || r | c | c   }
\hline
\textnumero  & Characteristic & Type  &\textnumero  & Characteristic & Type   \\  \hhline{=|=|= =|=|=}
$1$ &  {}$\ch010${} & $A_1$ & 			$8$ & {}$\ch222${}&  $A_3$         	 \\
$2$ &  {}$\ch101${} & $A_1$ & 			$9$ & {}$\ch121${} &  $C_2$         	 \\
$3$ &  {}$\ch200${} & $2A_1$ & 		$10$ & {}$\ch022${} &  $C_2+A_1$         	 \\
$4$ &  {}$\ch002${} & $2A_1$ & 		$11$ & {}$\ch224${} &  $C_2+A_1$         	 \\
$5$ &  {}$\ch111${} & $3A_1$ & 		$12$ & {}$\ch224${} &  $C_3$         	 \\
$6$ &  {}$\ch020${} & $2A_1$ & 		$13$ & {}$\ch422${} &  $C_3$         	 \\
$7$ &  {}$\ch202${} & $A_2$ & 		
\end{tabular}
\end{table}
\subsection{Mixed four-vectors, eighteenth family} Assume that the decomposition of the space $V$ into a direct sum of four two-dimensional subspaces is given: 
\[V=V_1 \oplus V_2 \oplus V_3 \oplus V_4.
\]

Then the space $V_1\otimes V_2\otimes V_3\otimes V_4$ can be identified as subspace of $\bw4 V$, assuming that for $v_k \in V_k$
\[v_1\otimes v_2\otimes v_3 \otimes v_4 = v_1 \wedge v_2 \wedge v_3 \wedge v_4.\]
The Lie algebra $L_0(V_1)\oplus L_0(V_2) \oplus L_0(V_3)\oplus L_0(V_4)$ can be naturally identified with a subalgebra in $L_0(V)$. The subspace 
\[\g^{(18)} =L_0(V_1)\oplus L_0(V_2) \oplus L_0(V_3)\oplus L_0(V_4) \oplus (V_1\otimes V_2\otimes V_3\otimes V_4)\]
 is a graded simple subalgebra of type $D_4$ in $\g$. 
 Suppose that 
$ \{e_1, e_2\} ,  \{e_3, e_4\} ,  \{e_5, e_6\} ,  \{e_7, e_8\} $ are bases of the subspaces $V_1, V_2, V_3, V_4$, respectively. In this case the algebra $\g^{(18)}$ is regular with respect to the Cartan subalgebra $\d< \g$, and $\g^{(18)}_0$ is a regular semisimple subalgebra of type $4A_1$. 

Let $p$ be the canonical semisimple four-vector of the 18th family, i.e., $p\in \c$ and $\c_p=\langle p_1, p_3, p_7\rangle$. It is easy to check that the simple component of the centralizer of the four-vector $p$ in the algebra coincides with the subalgebra $\g^{(18)}$ constructed above. 

Table~\ref{table:8} lists the classes of nilpotent four-vectors on the subalgebra $\g^{(18)}$. The subspace $\d^{(18)} = \d\cap  \g^{(18)}_0$ is taken as the Cartan subalgebra in $\g^{(18)}_0$, and the following  are chosen as simple roots of  $\g^{(18)}$: $\ep_1-\ep_2,\; \ep_3-\ep_4,\; \ep_5-\ep_6,\; \ep_7-\ep_8$.

\begin{table}
\caption{Nilpotent parts of four-vectors of the 18th family}\label{table:8}
\begin{tabular}{r | c | c || r | c | c   }
\hline
\textnumero  & Characteristic & Type  &\textnumero  & Characteristic & Type   \\  \hhline{=|=|= =|=|=}
$1$ &  {}$\ch1111${} & $A_1$ & 		$7$ & {}$\ch4000${}&  $A_2$         	 \\
$2$ &  {}$\ch2200${} & $2A_1$ & 		$8$ & {}$\ch4224${} &  $A_3$         	 \\
$3$ &  {}$\ch2020${} & $2A_1$ & 		$9$ & {}$\ch4242${} &  $A_3$         	 \\
$4$ &  {}$\ch2002${} & $2A_1$ & 		$10$ & {}$\ch4422${} &  $A_3$         	 \\
$5$ &  {}$\ch3111${} & $3A_1$ & 		$11$ & {}$\ch4844${} &  $D_4$         	 \\
$6$ &  {}$\ch2222${} & $4A_1$ & 		$12$ & {}$\ch0444${} &  $D_4(a_1)$         	 \\
\end{tabular}
\end{table}
\subsection{Mixed four-vectors, nineteenth family}\label{sec:5.8}
Let $p$ be the canonical semisimple four-vector of the 19th family, i.e., $p\in \c$ and $\c_p =  \langle p_1, p_1+p_3-p_7, p_2-p_3+p_6 \rangle$. The semisimple component $\g^{(19)}$ of the centralizer of the four-vector $p$ in algebra $\g$ (which is a simple subalgebra of type $A_4$) coincides with the intersection of the subalgebras $\g^{(2)}$ and $\g^{(12)}$ described respectively in Sections~\ref{sec:5.2} and \ref{sec:5.6} (see Table~\ref{table:1}), that is,
\[\g^{(19)}_0 = \{ X\in \g^{(12)}_0 \mid X^s_iB_{sj} + X^s_iB_{is} =0\},\]
\[\g^{(19)}_1 = \{ T\in \g^{(12)}_1 \mid T^{srij}B_{sr}=0\},\] 
with $B = e^1\wedge e^2 + e^3 \wedge e^4 + e^5\wedge e^6 + e^7 \wedge e^8$.

Let $\{f^1, f^2, f^3, f^4\}$ be the basis of the space $W^*$ dual to the basis $\{f_1, f_2, f_3, f_4\}$ of the space $W$. It is easy to check that 
\[\g^{(19)}_0 = \{A \in L_0(W) \mid A^s_k Q_{si} + A^s_i Q_{ks} =0  \},\]
\[\g^{(19)}_1 = \{T \in S^2(\bw2W) \mid T^{ijkl}d_{ijkl}=0,\;\; T^{ijkl}Q_{ij} =0  \},\]
with $d$ an arbitrary nonzero element of the space $\bw4W^*$, 
\[Q = f^1\wedge f^2 + f^3+f^4 \in \bw2W^*.\]
In particular, $\g^{(19)}_0$ is a simple algebra of type $C_2$.

Table~\ref{table:9} lists the classes of nilpotent four-vectors of the subalgebra $\g^{(19)}$.  The subspace $\d^{(19)} = \d \cap \g^{(19)}_0$ is taken as the Cartan subalgebra in $\g^{(19)}_0$, and the following are taken as simple roots with respect to $\d^{(19)}_0$: $\ep_1 -\ep_3, \ep_4-\ep_4$.

\begin{table}
\caption{Nilpotent parts of four-vectors of the 19th family}\label{table:9}
\begin{tabular}{r | c | c || r | c | c   }
\hline
\textnumero  & Characteristic & Type  &\textnumero  & Characteristic & Type   \\  \hhline{=|=|= =|=|=}
$1$ &  $0\;1$ & $A_1$ & 		$4$ & $1\;1$&  $2A_1$         	 \\
$2$ &  $1\;0$ & $A_1$ & 		$5$ & $1\;2$&  $C_2$         	 \\
$3$ &  $0\;2$ & $A_1$ & 		$6$ & $2\;2$&  $C_2$         	 \\
\end{tabular}
\end{table}

\section{Translators' Appendix I: Constructing normal forms of nilpotent vectors.}
In contrast to \cite{Vinberg-Elasvili}  Antonyan \cite{Antonyan} did not provide normal forms associated with each family of nilpotent vectors, rather he just lists the characteristic and the type and dimension of stabilizer subalgebra. Though the method for finding a ``canonical'' normal form for each nilpotent orbit is completely described in the article, we find it useful to just have a list.  It is a matter of computer programming and careful book-keeping. We automated this process in Macaulay2 \cite{M2}, and list the results in Table~\ref{table:nnForms}. 
We note that it is also possible to do this using de Graaf's \texttt{GAP} \cite{GAP4} packages \cites{SLA1.5.3, CoReLG1.56, QuaGroup1.8.3}, and whose method is described at \cite{de2011computing}, however we note that the representatives provided by this program are in terms of canonical bases and might not have the same form as we present here. Also, we note that our method essentially only uses linear algebra with some careful variable substitutions to reduce the number of free parameters in the under-determined systems that arise.
Below we explain what we did and provide some illustrative examples and connect the normal forms to the theory of Carter diagrams.

As described in Section~4.1\cite{Vinberg-Elasvili} a ``canonical'' element with dominant characteristic $h$ is the sum of root vectors corresponding to linearly independent roots, the number of which is the same as the rank of the algebra $\h$.

Let us review how one can find a nilpotent element associated with a given characteristic. 
We consider the families of nilpotent orbits in Table~\ref{table:2}. Each characteristic $h$ should coincide with a traceless diagonal matrix (being an element in $\d \subset L_0 V = \sl(V)$). Noting the fact that $\d$ is a Cartan subalgebra of a simple Lie algebra of type $A_7$, a set of $7$ simple roots of $\g$ with respect to $\d$ are chosen to be $\Pi_0 = \{\ep_1-\ep_2, \ldots, \ep_7-\ep_8\}$, however, given the fact that $\sum_i\ep_i=0$ we have a basis of $\d$ of elements of the form $d_i = \frac{1}{8} \diag((8-i)^i, (i)^{8-i})$ (here exponent denotes repetition), for $1\leq i\leq 7$.

The inner product on roots which is consistent with \eqref{eq:1.7} is 
\begin{equation}\label{eq:ip}
(
\ep_i + \ep_j + \ep_k + \ep_l, \ep_p + \ep_q + \ep_r + \ep_s) = d-2,
\end{equation}
 with $d$ the number of elements in common between the two sets $\{i,j,k,l\}$ and $\{p,q,r,s\}$. This means that for root vectors $e_i \wedge e_j \wedge e_k \wedge e_l$ and $e_p\wedge e_q \wedge e_r \wedge e_s$ the bracket from \eqref{eq:WWaction} reduces to 
 \[
 [e_i \wedge e_j \wedge e_k \wedge e_l, e_p\wedge e_q \wedge e_r \wedge e_s] 
 =(\ep_i + \ep_j + \ep_k + \ep_l, \ep_p + \ep_q + \ep_r + \ep_s) 
. \]
 This, together with the bilinearity of the bracket allows for quick computation of the brackets $[\;,\;] \colon \bw 4 \CC^8 \times \bw 4 \CC^8 \to \g$. Likewise, the formula at \eqref{eq:bracket01} allows for quick computation of brackets $[\;,\;] \colon \g \times \bw 4 \CC^8 \to \bw 4 \CC^8$. 

Following the support method, and the consequence of Morozov's theorem, given $h$, to find $e \in \bw{4} V$ we must solve the system of equations given by 
\begin{equation}\label{eq:Morozov}
[h, e] =2e,\quad  [h, f]=-2f,\quad  [e, f]= h
\end{equation}
 where the brackets are described at \eqref{eq:LWaction} and \eqref{eq:WWaction}. 
 The conditions $[h, e] =2e$, $[h, f]=-2f$ say that the pair $(e,f)$ is an element of  $E_{2} \times E_{-2}$, where $E_\lambda \subset \bw4 \CC^8$ is the $\lambda$-eigenspace of the adjoint operator associated with $h$.  
Thus the system at \eqref{eq:Morozov} reduces to solving  $[e, f]= h$  as a system of bilinear equations on $E_{2} \times E_{-2}$. Moreover, the computation of a basis of $E_\lambda$ for a given $h$ and its matrix representative $A$ is almost automatic given \eqref{eq:bracket01}; one just selects root vectors (standard basis vectors) $e_i \wedge e_j \wedge e_k \wedge e_l $ such that the partial trace (using diagonal entries given by the indices $(i,j,k,l)$)  of $A$ is equal to $\lambda$.

The bilinear system $[e,f]=h$ is often under determined and there can be many solutions.  Often this has a solution with $f = e^*$, the Hodge dual of $e$, but this is not always the case.  We can also seek a solution with support size equal to the rank of the stabilizer sub-algebra associated with that orbit, there exist solutions with rank equal to the rank of $\s$.
We report the results of these computations in Table~\ref{table:nnForms}.  We remark that Antonyan gives the list of nilpotent elements by a non-contiguous set of labels $1, \,\ldots, 94$. There are only 62 elements on the list, and Antonyan explains that he omits some orbits and explains his reasoning. Using de Graaf's \texttt{GAP} \cite{GAP4} packages \cites{SLA1.5.3, CoReLG1.56, QuaGroup1.8.3}. 
We noticed that when an orbit is omitted the Characteristic for orbit \textnumero $(i+1)$ should be the reverse of that for orbit \textnumero $i$.  However, by that reasoning orbit \textnumero 94 might also have been omitted because its character is the reverse of the character for orbit \textnumero 93.  One notices that in Antonyan's list of characteristics the characteristic is either palindromic or the subsequent orbit is omitted, with only one exception between orbits \textnumero 83 and \textnumero 86. Perhaps orbit \textnumero 85 (with characteristic $\ch440444 $) should have been retained and \textnumero 86 omitted. Because of this potential confusion we provide normal forms for all 94 orbits.

Let us explain further what happens when finding the normal forms for orbits \textnumero 93 and 94. Respectively the ``e'' and ``f'' for each are:

\[ \begin{matrix}e= e_{1}e_{2}e_{3}e_{4}+e_{0}e_{2}e_{4}e_{5}+e_{0}e_{3}e_{4}e_{6}+e_{0}e_{1}e_{5}e_{6}+e_{0}e_{2}e_{3}e_{7}+e_{0}e_{1}e_{4}e_{7},\\
f = \:e_{2}e_{3}e_{5}e_{6}+3\,e_{1}e_{4}e_{5}e_{6}+4\,e_{2}e_{3}e_{4}e_{7}-e_{1}e_{2}e_{5}e_{7}-e_{1}e_{3}e_{6}e_{7}+3\,e_{0}e_{5}e_{6}e_{7},
\end{matrix}
\]

\[ \begin{matrix}
\tilde e= e_{1}e_{2}e_{4}e_{5}+e_{0}e_{3}e_{4}e_{5}+e_{1}e_{2}e_{3}e_{6}+e_{0}e_{1}e_{4}e_{6}+e_{0}e_{2}e_{5}e_{6}+e_{0}e_{1}e_{2}e_{7},\\
\tilde f = \:3\,e_{3}e_{4}e_{5}e_{6}-e_{1}e_{3}e_{4}e_{7}-e_{2}e_{3}e_{5}e_{7}+3\,e_{0}e_{4}e_{5}e_{7}+4\,e_{1}e_{2}e_{6}e_{7}+e_{0}e_{3}e_{6}e_{7}.
\end{matrix}\]

We can move $e$ to $\tilde f$ by an element of $\GL_8$ as follows. Apply the substitution $e_i \mapsto e_{7-i}$ to $e$ to obtain a vector with the same support as $\tilde f$, and one can adjust the coefficients to match also by applying  the element $\diag \left\{-\frac{1}{4}t,\:1,\:-\frac{3}{4}t,\:1,\:-\frac{1}{4}t,\:t,\:-1,\:1\right\}$, with $t = \frac{ \pm 2}{\sqrt 3 }$. 

This computation and the evidence seen in Table~\ref{table:nnForms} and Figure~\ref{hasse} suggests the following.
\begin{conjecture} Using the notation in this article.
Suppose $e,f,h$ form an $\sl_2$ triple. 
\begin{enumerate}
\item The orbits of $e$ and $f$ are isomorphic.
\item If $\overleftarrow h$ has the reverse characteristic as $h$ with $\sl_2$ triple $\tilde e, \tilde f,\overleftarrow h$ has $\tilde f$ in the orbit of $e$. 
\item The characteristic of $h$ is palindromic if and only if $f$ is in the orbit of $e$. 
\end{enumerate}
\end{conjecture}

For some cases we found a different dimension of the orbit than what is in the original table of Antonyan, and we made the correction in the text in Table~\ref{table:2}. We summarize the discrepancies we noticed below. We did this computation using geometry, calculating the tangent space of the orbit, and later confirmed these results using de Graaf's \texttt{GAP} functions. 

 \begin{tabular}{c|c|c}
\textnumero  & dim  &  claimed dim  \\  \hhline{=|=|=}
      22&54&56\\
      36&46&49\\
      40&54&55\\
      43&45&48\\
      49&52&53\\
      72&49&52\\
\end{tabular}
 
Following are several examples which illustrate the computations we did. All our computations were done in Macaulay2 \cite{M2}.

We can place root vectors in a normal form on a Dynkin diagram with angles between vectors indicated by edges of different types. When not all angles are obtuse one calls the diagram a Carter diagram. 
We indicate by a dashed edge when the corresponding inner product of roots equal to $1$, no edge when the inner product is equal to 2, and a solid edge when the inner product is $-1$, (see \cite{Stekolshchik2010RootSA} for a nice exposition and many useful tables). Note: The notation $B(0,1,2)$ comes from the language of Shephard and Todd \cite{ShephardTodd54}.

\begin{example} Consider the characteristic $h = \ch0001000$ for the nilpotent family \textnumero $1$ Table~\ref{table:2}. 
In this case $h =\diag(\frac{1}{2},\frac{1}{2},\frac{1}{2},\frac{1}{2}, \frac{-1}{2}, \frac{-1}{2}, \frac{-1}{2}, \frac{-1}{2}) $.

We have 
\[
[h,e_i \wedge e_j \wedge e_k \wedge e_l] = \sum_{1\leq i<j<k<l\leq 8} (\ep_i + \ep_j + \ep_k + \ep_l)(h) e_i \wedge e_j \wedge e_k \wedge e_l
.\]
In this case
\[
(\ep_i + \ep_j + \ep_k + \ep_l)(h)
= \frac{1}{2} \cdot (\#(\{1,2,3,4\}\cap \{i,j,k,l\} ) - \#(\{5,6,7,8\}\cap \{i,j,k,l\} ) 
.\]
We see that  $e = e_1\wedge e_2 \wedge e_3 \wedge e_4$ is the only eigenvector with eigenvalue $2$, and $f = e^* = e_5 \wedge e_6 \wedge e_7 \wedge e_8$, and indeed $[e,e^*] = h$. 
\end{example}

\begin{example}
Consider the characteristic $h = \ch0102010$ for the nilpotent family \textnumero $20$ Table~\ref{table:2}. 
In this case $h =\diag(2,2,1,1,-1,-1,-2,-2) $.

In this case $(\ep_i + \ep_j + \ep_k + \ep_l)(h)
=
$
\[
 2\cdot \#(\{1,2\}\cap \{i,j,k,l\} ) 
+\#(\{3,4\}\cap \{i,j,k,l\} ) 
- \#\{5,6\}\cap \{i,j,k,l\} )
-  2\cdot \#\{7,8\}\cap \{i,j,k,l\} ) 
.\]
We see that  the vectors  $\left\{e_{1}e_{2}e_{5}e_{6},\,e_{1}e_{3}e_{4}e_{7},\,e_{2}e_{3}e_{4}e_{7},\,e_{1}e_{3}e_{4}e_{8},\,e_{2}e_{3}e_{4}e_{8}\right\}$ are the basis elements with eigenvalue $2$. So any linear combination of them will have eigenvalue 2 for $\ad(h)$.  The $-2$ eigenspace is spanned by $\{e_1e_5e_6e_7,\, e_2e_5e_6e_7,\, e_1e_5e_6e_8,\, e_2e_5e_6e_8,\, e_3e_4e_7e_8\}$.
By standard polynomial algebra one can find a restriction of the variety of solutions to $[e,f]=h$ that has rank-3. One such is 
\[
e_{1256}+e_{1347}+e_{2348},
\]
where we use the standard shorthand $e_{ijkl} = e_ie_je_ke_l$. 
The terms in this form can be arranged as labels for a Dynkin diagram for $A_3$:
\[
\xymatrix{
\ccircle{1347} \ar@{-}[r] & \ccircle{1256}\ar@{-}[r] &\ccircle{2348} 
}
\]
\end{example}
\begin{example}
Consider the characteristic $h = \ch1211121$ for the nilpotent family \textnumero $22$ Table~\ref{table:2}. 
In this case $h =\diag(2,2,1,1,-1,-1,-2,-2) $, which corresponds to the diagonal matrix $\diag  \left( \frac{9}{2},\,\frac{7}{2},\,\frac{3}{2},\,\frac{1}{2},\,-\frac{1}{2},\,-\frac{3}{2},\,-\frac{7}{2},\,-\frac{9}{2}\right)$

In this case there are precisely $5$ basis elements with eigenvalue $2$, and 5 with eigenvalue $-2$ for $h$. We find that 
\[\begin{matrix}
e &=& e_{2456}+e_{2347}+e_{1357}+e_{1348}+e_{1268},\\
f &=& -8e_{3457}+5e_{2567}+9e_{2468}+5e_{1568}+-8e_{1378},
\end{matrix}
\]
solve $[e,f] = h$.
The terms in this form can indeed be arranged as labels for a Dynkin diagram for $A_5$:
\[
\xymatrix{
\ccircle{1348} \ar@{-}[r] & \ccircle{2456} \ar@{-}[r] & \ccircle{1357} \ar@{-}[r] & \ccircle{1268}  \ar@{-}[r] & \ccircle{2347}
}
\]
We have a parametrization of the orbit of  a normal form $x$ via $G \times {x} \to \bw{4} \CC^8$ given by $(g,x) = g.x$, and we can use this to compute the dimension of the orbit as the rank of the Jacobian at the point $(e.x)$.  When we compute the dimension of the tangent space at a general point of the orbit, however, we find that this orbit has dimension $54$ and not $56$ as would be indicated by the dimension of the stabilizer subalgebra in table \ref{table:2}.
\end{example}

\begin{example}
For nilpotent orbit \textnumero $88$, we find normal form:
\[e_{3456}+e_{2457}+e_{2367}+e_{1467}+e_{2348}+e_{1358}+e_{1268}.\]
The terms in this form can indeed be arranged as labels for the Dynkin diagram for $E_7$ 
\begin{equation*}\belowdisplayskip=-12pt
\xymatrix{
&&&\ccircle{2367}\ar@{-}[d]&&&
\\
\ccircle{3456} \ar@{-}[r] & \ccircle{1268}\ar@{-}[r]  &\ccircle{2457}\ar@{-}[r]  &\ccircle{1358} \ar@{-}[r]  &\ccircle{1467} \ar@{-}[r]  &\ccircle{2348} 
}
\end{equation*}
\qedhere
\end{example}
\begin{example}
To interpret types indicated by $G(\alpha)$ we refer to \cite[Theorem~5.5]{dynkin1960semisimple} and the discussion surrounding it. According to Dynkin for a semi-simple group $G$ with root system $\Sigma$ and simple roots $\Pi$ we let $\Pi_\a$ denote the system obtained by omitting $\a$, and $G(\alpha)$ denotes the subalgebra of $G$ generated by $\{e_\b, e_{-\b} \mid \b \in \Pi_\a\} \cup \{e_\delta, e_{-\delta}\}$ with $\delta$ the smallest root of $G$ not expressible in terms of $\pm \b$ with $\b \in \Pi_\a$. 

For example, for orbit  \textnumero 44, which has type $D_4(a_1)$, we find normal form:
\[
e_{1246}+e_{1357}+e_{2348}+e_{2358}
.\]
We start with a labeling of a diagram of type $A_3$ using roots in the $2$-eigenspace associated with characteristic $h = \ch0020200$, as:
\[
\xymatrix{
&&\ccircle{2358}\\
& \ccircle{1246} \ar@{-}[ur]  \ar@{-}[dr]\\
&& \ccircle{1357} 
}
\]
We can choose $\delta = \ep_4-\ep_5$ and append the node associated with weight $(\ep_2 +\ep_3 +\ep_5 + \ep_8) +(\ep_4-\ep_5)= (\ep_2 +\ep_3 +\ep_4 + \ep_8) $ to obtain the Carter diagram:
\[
\xymatrix{
&&\ccircle{2358}\\
& \ccircle{1246} \ar@{-}[ur]  \ar@{-}[dr] &&
\ccircle{2348} \ar@{--}[lu] \ar@{-}[dl] 
\\
&& \ccircle{1357} 
}
\]
which indeed follows the pattern for $D_4(a_1)$ seen in \cite[Figure~1.1]{Stekolshchik2010RootSA}.
We also notice that the normal form has a simpler (from the point of view of tensor rank) format since
\[e_1e_2e_4e_6+e_1e_3e_5e_7+e_2e_3e_4e_8+e_2e_3e_5e_8
=
e_1e_2e_4e_6+e_1e_3e_5e_7+e_2e_3(e_4+e_5)e_8,
\]
which we might represent as the following diagram (because the inner product on roots should be bilinear):
\[
\xymatrix{
&&\ccircle{23(4+5)8}\ar@{-}[dd]\\
& \ccircle{1246} \ar@{-}[ur]  \ar@{-}[dr] &&
\\
&& \ccircle{1357} 
}
\]
From this we might ask if it makes sense to collapse all the acute angles in the Carter diagrams for the other nilpotent orbits of type $G(a)$, and thus obtain a different type of diagram. 
\end{example}
\begin{example} Consider the normal form for \textnumero42,
$e_{1456}+e_{2347}+e_{1267}+e_{1358}+e_{2358}+e_{1268}$, which is supposed to have type $E_6(b)$. 
We might arrange the vectors in the Carter diagram:
\[
\xymatrix{
&&\ccircle{1268} \ar@{--}[rd] \ar@{-}[ld] \\
\ccircle{1456} \ar@{-}[r] \ar@{-}@/_2pc/[rrrr] & \ccircle{2347}\ar@{-}[r]  &\ccircle{1358}\ar@{-}[r]\ar@{--}@/_2pc/[rr]  &\ccircle{1267} \ar@{-}[r]  &\ccircle{2358}
\\
\\}
\]
\bigskip
This diagram has partial Cartan matrix 
$\left(\begin{smallmatrix}
      2&-1&0&0&-1&0\\
      -1&2&0&-1&0&-1\\
      0&0&2&-1&-1&1\\
      0&-1&-1&2&1&0\\
      -1&0&-1&1&2&0\\
      0&-1&1&0&0&2\\
      \end{smallmatrix}\right)$.

We could contract the dashed arrows to produce the normal form
$e_1e_4e_5e_6+e_2e_3e_4e_7+e_1e_2e_6(e_7+e_8)+(e_1+e_2)e_3e_5e_8$, and corresponding diagram:
\[
\xymatrix{
 \ccircle{126(7+8)} \ar@{->}[d] \ar@{=}[r]&\ccircle{(1+2)358}\ar@{->}[dl]\ar@{->}[d]\\
\ccircle{2347}  \ar@{-}[r] &\ccircle{1456} 
}
\]
The arrows point from the longer to the shorter roots when they are not the same length.
The contracted Cartan matrix is:
$\left(\begin{smallmatrix}
      2&-1&0&-1\\
      -1&2&-1&-1\\
      0&-1&6&-2\\
      -1&-1&-2&6\\
      \end{smallmatrix}\right)$.
\end{example}

\subsection{Orbit Closures of Nilpotent Elements}
At the suggestion of a referee, we performed some additional computations utilizing the \texttt{GAP} \cite{GAP4} packages \cites{SLA1.5.3, CoReLG1.56, QuaGroup1.8.3}. We found \cite{deGraafclosure} as well as \cites{Popov03, de2017computation,DerksenKemper }  also to be helpful in understanding these calculations.
 The Hasse diagram for Antonyan's labels as calculated by de Graaf's program and converted to Antonyan's indexing consists of the directed edges listed in Figure~\ref{edges}. We depict the Hasse diagram in Figure \ref{hasse}.

%
%

\section*{Acknowledgements}
I would like to thank  Fr\'ed\'eric Holweck for useful discussions. Antonyan's original paper was difficult to find, and I am thankful for Professor Vladimir Popov for supplying me a copy. I thank Willem de Graaf for help in utilizing the code he wrote.

I also note that I am not competent in Russian, however, as my collaborator and I found ourselves referring to this paper very often for recent work, we decided that it would be useful to us, and hopefully to the greater community of mathematicians and physicists, to have an English translation.  I am grateful to Google Translate for help with translating the prose. Unfortunately at the time of this writing it was not possible to just upload the scanned copy of the article to Google Translate's website and obtain a completely translated document, instead it was possible to use a phone camera and Google Lens and holding a mobile device up to the computer screen. Obviously this is not optimal, and needs not just human input, but also subject knowledge, and this is one reason to try to translate the article. I have done my best to fill in the gaps left by the automatic translator. I checked the mathematical meaning of the translated prose based on the fact that the methods follow very closely the work of Vinberg and company \cites{GV, AntonyanElashvili}, which is well-explained in English in \cite{Vinberg-Elasvili}. In essence, this translation is a copy of the mathematical formulas from Antonyan, and the methods from Vinberg, and the prose from Google Translate.

\begin{table}
\caption{Canonical Nilpotent forms in $\bw{4}\CC^8$}\label{table:nnForms}
\resizebox{.9 \textwidth}{!}{
\begin{tabular}{c | c | l | r}
\textnumero  & characteristic &  normal form  & dim   \\  \hhline{=|=|=|=}
1& \ch0001000 &$e_{1234}$&17\\
2& \ch0100010 &$e_{1234}+e_{1256}$&26\\
3& \ch0200000 &$e_{1234}+e_{1256}+e_{1278}$&27\\
4& \ch0000020 &$e_{1234}+e_{1256}+e_{3456}$&27\\
5& \ch1001001 &$e_{1235}+e_{1246}+e_{1347}$&32\\
6& \ch0002000 &$e_{1235}+e_{1246}+e_{1347}+e_{2348}$&33\\
7& \ch1100100 &$e_{1345}+e_{1236}+e_{1247}+e_{1258}$&35\\
8& \ch0010011 &$e_{1245}+e_{1346}+e_{2356}+e_{1237}$&35\\
9& \ch2000002 &$e_{1234}+e_{1567}$&33\\
10& \ch2010001 &$e_{1245}+e_{1367}+e_{1238}$&38\\
11& \ch1000102 &$e_{2345}+e_{1236}+e_{1457}$&38\\
12& \ch0101010 &$e_{1345}+e_{2346}+e_{1256}+e_{1237}+e_{1248}$&38\\
13& \ch3000100 &$e_{1246}+e_{1357}+e_{1238}+e_{1458}$&41\\
14& \ch0010003 &$e_{1245}+e_{1346}+e_{1357}+e_{2367}$&41\\
15& \ch1010101 &$e_{2345}+e_{1246}+e_{1357}+e_{1238}$&41\\
16& \ch4000000 &$e_{1246}+e_{1357}+e_{1238}+e_{1458}+e_{1678}$&42\\
17& \ch0000004 &$e_{1234}+e_{1256}+e_{1357}+e_{2367}+e_{4567}$&42\\
18& \ch2000200 &$e_{2345}+e_{1246}+e_{1357}+e_{1238}+e_{1458}$&42\\
19& \ch0020002 &$e_{1245}+e_{1346}+e_{1357}+e_{2367}+e_{1238}$&42\\
20& \ch0102010 &$e_{1256}+e_{1347}+e_{2348}$&42\\
21& \ch2200022 &$e_{2345}+e_{1347}+e_{1567}+e_{1268}$&50\\
22& \ch1211121 &$e_{2456}+e_{2347}+e_{1357}+e_{1348}+e_{1268}$&54\\
23& \ch0202040 &$e_{3456}+e_{1347}+e_{1257}+e_{2348}+e_{1268}$&51\\
24& \ch0402020 &$e_{1356}+e_{2456}+e_{1347}+e_{2348}+e_{1278}$&51\\
25& \ch2220222 &$e_{2456}+e_{2347}+e_{1457}+e_{1367}+e_{1358}+e_{1268}$&57\\
26& \ch2222222 &$e_{2456}+e_{2357}+e_{1457}+e_{1367}+e_{2348}+e_{1358}+e_{1268}$&59\\
27& \ch0200020 &$e_{1345}+e_{2346}+e_{1257}+e_{1268}$&42\\
28& \ch1030010 &$e_{2345}+e_{1347}+e_{1267}+e_{1258}+e_{1368}$&49\\
29& \ch0100301 &$e_{2346}+e_{1356}+e_{1347}+e_{2457}+e_{1258}$&49\\
30& \ch0040040 &$e_{3456}+e_{1247}+e_{1357}+e_{2358}+e_{1268}+e_{1368}+e_{2368}$&56\\
31& \ch0400400 &$e_{1346}+e_{2357}+e_{1267}+e_{1458}+e_{2458}+e_{1268}+e_{1278}$&56\\
32& \ch2020202 &$e_{2356}+e_{2347}+e_{1457}+e_{1367}+e_{1248}+e_{1358}$&53\\
33& \ch0202020 &$e_{1356}+e_{2456}+e_{1347}+e_{1257}+e_{2348}+e_{1268}$&50\\
34& \ch0002020 &$e_{1256}+e_{3456}+e_{1347}+e_{2348}$&43\\
35& \ch0202000 &$e_{1256}+e_{1347}+e_{2348}+e_{1278}$&43\\
36& \ch1011101 &$e_{1256}+e_{2347}+e_{1357}+e_{1348}$&46\\
37& \ch3101021 &$e_{2345}+e_{1357}+e_{1467}+e_{1348}+e_{1268}$&52\\
38& \ch1201013 &$e_{2356}+e_{2347}+e_{1357}+e_{1467}+e_{1248}$&52\\
39& \ch2202022 &$e_{2456}+e_{2347}+e_{1357}+e_{1467}+e_{1348}+e_{1268}$&55\\
40& \ch1311111 &$e_{2356}+e_{1456}+e_{2347}+e_{1357}+e_{1348}+e_{1278}$&54\\
\end{tabular}}
\end{table}
\begin{table}
\resizebox{.9 \textwidth}{!}{
\begin{tabular}{c | c | l | r}
\textnumero  & characteristic &  normal form  & dim   \\  \hhline{=|=|=|=}
41& \ch1111131 &$e_{3456}+e_{2347}+e_{1357}+e_{1267}+e_{1348}+e_{1258}$&54\\
42& \ch0220220 &$e_{1456}+e_{2347}+e_{1267}+e_{1358}+e_{2358}+e_{1268}$&55\\
43& \ch1101011 &$e_{2345}+e_{1356}+e_{1347}+e_{1267}+e_{1248}$&45\\
44& \ch0020200 &$e_{1246}+e_{1357}+e_{2348}+e_{2358}$&47\\
45& \ch2020020 &$e_{2345}+e_{1456}+e_{1347}+e_{1267}+e_{1258}+e_{1368}$&50\\
46& \ch0200202 &$e_{2346}+e_{1356}+e_{1347}+e_{2457}+e_{1267}+e_{1258}$&50\\
47& \ch0040000 &$e_{2345}+e_{1347}+e_{1267}+e_{1258}+e_{1368}+e_{2378}$&50\\
48& \ch0000400 &$e_{1236}+e_{1456}+e_{1347}+e_{2457}+e_{1258}+e_{3458}$&50\\
49& \ch1111111 &$e_{2356}+e_{1456}+e_{2347}+e_{1357}+e_{1267}+e_{1348}+e_{1258}$&52\\
50& \ch2004002 &$e_{1256}+e_{1357}+e_{1467}+e_{2348}$&48\\
51& \ch2204022 &$e_{2456}+e_{1357}+e_{1467}+e_{2348}+e_{1268}$&56\\
52& \ch1313143 &$e_{3456}+e_{2357}+e_{1457}+e_{1367}+e_{2348}+e_{1268}$&59\\
53& \ch3413131 &$e_{2456}+e_{1457}+e_{1367}+e_{2348}+e_{1358}+e_{1278}$&59\\
54& \ch0202202 &$e_{1356}+e_{2457}+e_{1267}+e_{1348}+e_{2348}+e_{1258}$&54\\
55& \ch2022020 &$e_{2356}+e_{1456}+e_{2347}+e_{1367}+e_{1258}+e_{1368}$&54\\
56& \ch0004004 &$e_{1356}+e_{2457}+e_{1267}+e_{3467}+e_{1348}+e_{2348}$&54\\
57& \ch4004000 &$e_{2346}+e_{1257}+e_{1367}+e_{1358}+e_{1468}+e_{1478}$&54\\
58& \ch1013012 &$e_{2356}+e_{1456}+e_{1257}+e_{1367}+e_{2348}$&51\\
59& \ch2103101 &$e_{1356}+e_{1457}+e_{1267}+e_{2348}+e_{1258}$&51\\
60& \ch1112111 &$e_{2356}+e_{1456}+e_{1357}+e_{1267}+e_{2348}+e_{1258}$&53\\
61& \ch0103103 &$e_{1356}+e_{2457}+e_{1267}+e_{1348}+e_{2348}$&53\\
62& \ch3013010 &$e_{1456}+e_{2347}+e_{1367}+e_{1258}+e_{1368}$&53\\
63& \ch3113121 &$e_{2356}+e_{1457}+e_{1367}+e_{2348}+e_{1358}+e_{1268}$&57\\
64& \ch1213113 &$e_{2456}+e_{2357}+e_{1457}+e_{1367}+e_{2348}+e_{1258}$&57\\
65& \ch0404044 &$e_{3456}+e_{1357}+e_{2467}+e_{1348}+e_{2348}+e_{1258}+e_{1268}$&60\\
66& \ch4404040 &$e_{2356}+e_{2347}+e_{1457}+e_{1458}+e_{1368}+e_{1468}+e_{1278}$&60\\
67& \ch2422222 &$e_{2456}+e_{2357}+e_{1457}+e_{1367}+e_{2348}+e_{1358}+e_{1278}$&60\\
68& \ch2222242 &$e_{3456}+e_{2357}+e_{1457}+e_{1367}+e_{2348}+e_{1358}+e_{1268}$&60\\
69& \ch4224224 &$e_{2456}+e_{2357}+e_{1467}+e_{2348}+e_{1358}+e_{1268}$&60\\
70& \ch3013131 &$e_{2456}+e_{1457}+e_{1367}+e_{2348}+e_{1258}+e_{1358}$&57\\
71& \ch1313103 &$e_{2356}+e_{1457}+e_{1367}+e_{2348}+e_{1268}+e_{1278}$&57\\
72& \ch1110111 &$e_{2346}+e_{1456}+e_{1357}+e_{1267}+e_{1248}+e_{1258}$&49\\
73& \ch1010210 &$e_{2346}+e_{1456}+e_{1257}+e_{1348}+e_{1358}$&48\\
74& \ch0120101 &$e_{1346}+e_{2357}+e_{1267}+e_{1248}+e_{1258}$&48\\
75& \ch2022222 &$e_{2456}+e_{2357}+e_{1457}+e_{1367}+e_{2348}+e_{1258}+e_{1358}$&58\\
76& \ch2222202 &$e_{2356}+e_{1457}+e_{1367}+e_{2348}+e_{1358}+e_{1268}+e_{1278}$&58\\
77& \ch4004040 &$e_{2456}+e_{1357}+e_{1267}+e_{2348}+e_{1458}+e_{1368}+e_{1468}$&58\\
78& \ch0404004 &$e_{2356}+e_{1457}+e_{1367}+e_{2467}+e_{2348}+e_{1268}+e_{1278}$&58\\
79& \ch4220224 &$e_{2456}+e_{2347}+e_{1567}+e_{1348}+e_{1358}+e_{1268}$&59\\
80& \ch1030131 &$e_{2456}+e_{2347}+e_{1367}+e_{1348}+e_{1258}+e_{1358}$&55\\
\end{tabular}}
\end{table}
\begin{table}
\resizebox{.9 \textwidth}{!}{
\begin{tabular}{c | c | l | r}
\textnumero  & characteristic &  normal form  & dim   \\  \hhline{=|=|=|=}
81& \ch1310301 &$e_{2346}+e_{2357}+e_{1457}+e_{1358}+e_{1268}+e_{1278}$&55\\
82& \ch2020222 &$e_{2456}+e_{2347}+e_{1457}+e_{1367}+e_{1348}+e_{1258}+e_{1358}$&56\\
83& \ch4444044 &$e_{2456}+e_{2357}+e_{1467}+e_{2348}+e_{1358}+e_{1368}+e_{1278}$&62\\
84& \ch4404444 &$e_{3456}+e_{2357}+e_{1467}+e_{2348}+e_{1358}+e_{1458}+e_{1268}$&62\\
85& \ch2220202 &$e_{2346}+e_{2357}+e_{1457}+e_{1367}+e_{1358}+e_{1268}+e_{1278}$&56\\
86& \ch4404404 &$e_{2356}+e_{2457}+e_{1367}+e_{2348}+e_{1458}+e_{1268}+e_{1278}$&61\\
87& \ch4044044 &$e_{2456}+e_{2357}+e_{1467}+e_{2348}+e_{1358}+e_{1268}+e_{1368}$&61\\
88& \ch4444448 &$e_{3456}+e_{2457}+e_{2367}+e_{1467}+e_{2348}+e_{1358}+e_{1268}$&63\\
89& \ch8444444 &$e_{2456}+e_{2357}+e_{1467}+e_{2348}+e_{1458}+e_{1368}+e_{1278}$&63\\
90& \ch0101111 &$e_{1356}+e_{2456}+e_{2347}+e_{1257}+e_{1248}$&47\\
91& \ch1111010 &$e_{2346}+e_{1356}+e_{1257}+e_{1348}+e_{1268}$&47\\
92& \ch2002002 &$e_{1256}+e_{2347}+e_{1357}+e_{1467}+e_{1348}$&47\\
93& \ch2101101 &$e_{2345}+e_{1356}+e_{1457}+e_{1267}+e_{1348}+e_{1258}$&48\\
94& \ch1011012 &$e_{2356}+e_{1456}+e_{2347}+e_{1257}+e_{1367}+e_{1238}$&48\\
\end{tabular}
}
\end{table}
\begin{figure}[b]
 \caption{Edges in the Hasse Diagram of the orbit closures of nilpotent elements in $\bw{4}\CC^{8}$.}\label{edges}
((83, 89), (83, 88), (84, 89), (84, 88), (86, 83), (86, 84), (87, 83), (87, 84), (69, 86), (69, 87), (66, 86), (67, 86), (67, 87), (68, 86), (68, 87), (65, 87), (53, 66), (53, 67), (79, 69), (79, 66), (79, 65), (52, 68), (52, 65), (26, 69), (26, 67), (26, 68), (77, 53), (77, 79), (77, 52), (77, 26), (78, 53), (78, 79), (78, 52), (78, 26), (76, 53), (76, 26), (75, 52), (75, 26), (25, 77), (25, 78), (71, 78), (71, 76), (70, 77), (70, 75), (64, 78), (64, 75), (63, 77), (63, 76), (51, 71), (51, 70), (51, 64), (51, 63), (85, 25), (85, 71), (85, 63), (82, 25), (82, 70), (82, 64), (30, 25), (30, 71), (30, 70), (30, 64), (31, 25), (31, 71), (31, 70), (31, 63), (39, 51), (39, 85), (39, 82), (81, 85), (81, 31), (42, 51), (42, 30), (42, 31), (80, 82), (80, 30), (57, 85), (56, 82), (55, 81), (55, 42), (22, 39), (22, 81), (22, 42), (22, 80), (41, 80), (40, 81), (54, 42), (54, 80), (60, 55), (60, 54), (61, 56), (61, 54), (62, 57), (62, 55), (32, 55), (32, 22), (32, 41), (32, 40), (32, 54), (49, 60), (49, 32), (37, 62), (37, 32), (38, 61), (38, 32), (23, 41), (58, 60), (58, 61), (59, 60), (59, 62), (24, 40), (21, 37), (21, 38), (21, 23), (21, 24), (47, 49), (33, 49), (33, 23), (33, 24), (48, 49), (46, 49), (46, 38), (46, 58), (45, 49), (45, 37), (45, 59), (28, 47), (28, 33), (28, 45), (29, 33), (29, 48), (29, 46), (72, 21), (72, 33), (72, 46), (72, 45), (94, 46), (73, 29), (73, 72), (93, 45), (50, 58), (50, 59), (74, 28), (74, 72), (91, 93), (91, 74), (92, 72), (92, 94), (92, 93), (92, 50), (44, 73), (44, 50), (44, 74), (90, 94), (90, 73), (36, 91), (36, 92), (36, 44), (36, 90), (43, 36), (35, 91), (34, 90), (16, 93), (18, 43), (17, 94), (20, 36), (20, 35), (20, 34), (27, 43), (27, 35), (27, 34), (19, 43), (13, 16), (13, 18), (15, 18), (15, 20), (15, 27), (15, 19), (14, 17), (14, 19), (11, 15), (11, 14), (10, 13), (10, 15), (12, 15), (7, 10), (7, 12), (8, 11), (8, 12), (9, 11), (9, 10), (6, 12), (5, 7), (5, 8), (5, 9), (5, 6), (3, 7), (4, 8), (2, 5), (2, 3), (2, 4), (1, 2))
\end{figure}

%
\clearpage
\begin{figure}
\caption{The Hasse Diagram of the orbit closures of nilpotent elements in $\bw{4}\CC^{8}$ with heights corresponding to dimension. Orbits with palindromic characters are placed near the center, and non-palindromic characters have their reverse reflected over the horizontal center line.}
\label{hasse}
\includegraphics[height=0.9\textheight]{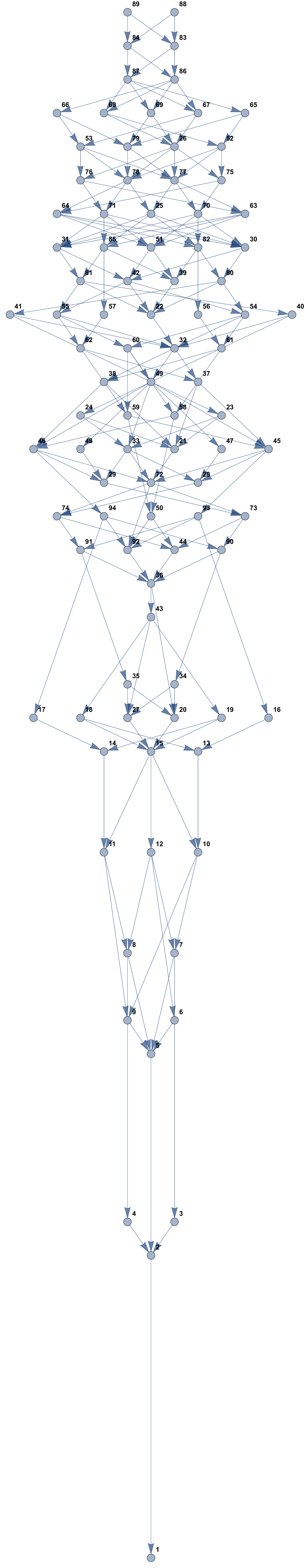} 
\end{figure}

 \clearpage
\bibliographystyle{amsplain}
\bibliography{/Users/oeding/Dropbox/BibTeX_bib_files/main_bibfile.bib}

\end{document}